\documentclass[11pt,a4paper,twoside]{amsart}
\usepackage{graphicx,calc,amscd}
\usepackage{epsfig,psfig}
\usepackage{float}
\usepackage{labelfig}
\usepackage{amsfonts}
\usepackage{amssymb}

\usepackage{fontenc}
\usepackage{amsmath}
\usepackage{amssymb}
\usepackage{latexsym}

\newcommand{\co}{\colon\thinspace}    
\newcommand{\fnote}[1]{\footnote{\small sharp1}}
\newcommand{\inv}{^{-1}}              

\newcommand{\N}{{\mathbb N}}
\newcommand{\K}{{\mathbb K}}
\newcommand{\Z}{{\mathbb Z}}
\newcommand{\R}{{\mathbb R}}

\newcommand{\T}{{\mathbb T}}

\newcommand{\spt}{\mbox{supp}}

\newcommand{\s}{\mathcal{S}_1}

\newcommand{\eps}{\bar{\epsilon}}
\newcommand{\inter}{\mbox{Int}}
\newcommand{\supp}{\mbox{supp}}
\newcommand{\vect}{\mbox{Vect}}

\newtheorem{theorem}{Theorem}
\newtheorem{proposition}[theorem]{Proposition}
\newtheorem{corollary}[theorem]{Corollary}

\newtheorem{lemma}[theorem]{Lemma}

\title{Stable norms of non-orientable surfaces}
\author{Florent Balacheff, Daniel Massart}
\date{\today}
\begin{document}

\begin{abstract}
We study the stable norm on the first homology of a closed, non-orientable surface equipped with a Riemannian metric. We prove that in every conformal class  there exists a metric whose stable norm is polyhedral. Furthermore the stable norm is never strictly convex if the first Betti number of the surface is greater than two.
\end{abstract}

\maketitle
\section{Introduction}
Given a compact Riemannian manifold $(M,g)$ with first Betti number $b_1(M) >0$, the stable norm $\|\ \|$ on $H_1(M,\R)$ is defined in \cite{GLP} (see also \cite{Federer}) as
    \[\begin{array}{rcl}
 H_1(M,\R) & \longrightarrow  & \R\\
h & \longmapsto & \|h\| := \inf \left\{\sum^{n}_{i=1}|r_i|l_g (\gamma_i)     \right\}
\end{array}
\]
where
\begin{itemize}
    \item
    $l_g$ denotes the length with respect to $g$
    \item
    the $r_i$ are real numbers
    \item
    the $\gamma_i$ are Lipschitz 1-cycles
    \item
    $h = \sum^{n}_{i=1} r_i \left[ \gamma_i\right]$.
\end{itemize}
Note that since we want to minimize the length we may assume from the start that the $\gamma_i$ are  closed geodesics that minimize the length in their free homotopy class.

In general the infimum may not be reached. It is remarkable that when the dimension of $M$ is two, it is reached for every integer homology
class (Proposition \ref{rational}). When the infimum is actually a minimum, we may wonder whether the minimizing cycles are connected. Note that
every component $\gamma_i,i=1\ldots n$ of a minimizing cycle $\sum^{n}_{i=1} r_i \gamma_i$ is itself minimizing in its own homology class. A
minimizing cycle whose connected components have distinct homology classes yields a flat region in the unit sphere $\s$ of the stable norm,
containing the convex hull of the $\{[\gamma_i] /l_g(\gamma_i)\}_{i=1}^n$. So we may ask how often does it occur, how many components do the
minimizing cycles have and what is the dimension of the corresponding flat (that is the dimension of the affine subspace it spans). In this
paper we give some answers when $M$ is a closed non-orientable surface. Our first result is similar to Theorem 7 of \cite{gafa} which adresses
the orientable case. We denote by $[x]$ the integer part of a real number $x$.

\bigskip
\noindent {\bf Theorem A} {\it Assume $M$ is a closed non-orientable surface  endowed with a Riemannian metric. Then every connected minimizing
cycle is a component of a minimizing cycle with  at least $[(b_1(M)+1)/2]-1$ homologically independant components, and at most $2b_1(M)-1$ non
pairwise homologically proportional components.}

\bigskip

So the stable norm is never strictly convex for $b_1(M)>2$. The difference with Theorem 7 of \cite{gafa} is that the dimension of the corresponding flat may be more than $[(b_1(M)+1)/2]-1$. Observe that if $b_1(M)=2$ the stable norm may be strictly convex. For instance, take a hyperbolic punctured torus, cut off a sufficiently thin neighborhood of the cusp, and glue a projective plane.

Let $\pi \co M_o \longrightarrow M$ be the orientation cover of a non-orientable surface $M$. A simple closed curve $\gamma$ of $M$ is said {\it of type I} (resp. {\it of type II}) if its inverse image $\pi^{-1}(\gamma)$ consists of either one curve or two homologous curves (resp. two non-homologous curves). Remark that one-sided simple closed curves on $M$ (curves whose tubular neighbourhood is homeomorphic to a M\H{o}bius strip) are of type I while two-sided simple closed curves on $M$ (curves whose tubular neighbourhood is homeomorphic to an annulus) may be of type I or II. The following theorem states that the local geometry of the unit sphere $\s$ is special near homology classes whose minimizing cycles consist of curves of type I. Specifically, the intersection of the unit ball with a neighborhood of such a class is a cone .

\bigskip

\noindent {\bf Theorem B}
{\it Assume $M$ is a closed non-orientable surface endowed with a Riemannian metric. Let $h_0$ be an integer homology class all of whose minimizing cycles consist of geodesics of type I. Then for all $h \in H_1(M,\R)$, there exists $s(h_0,h)>0$ such that the subset of the unit sphere $\s$
$$
 \left\{\frac{ h_0+ sh}{|| h_0+s h||} \co s \in \left[0,s(h_0,h)\right] \right\}
$$
is a straight segment.
}

\bigskip

Apart from surfaces little is known about minimizing cycles. For flat tori they exist and are connected in every integer homology class (or multiple thereof). Other homology classes do not have minimizing cycles.
Furthermore the stable norm of a flat torus is Euclidean. Apart from \cite{McShane} which deals with hyperbolic metrics on a punctured torus, the only other examples (\cite{Bangert90}, \cite{Babenko-Balacheff06}) where the stable norm is actually  computed have very few connected minimizing cycles : the unit ball of the stable norm is a polyhedron. So there is but a finite number of connected minimizing cycles, corresponding to the vertices of the polyhedron. In  every homology class there is a minimizing cycle which is a linear combination of the connected ones. All such examples assume $\dim M \geq 3$ ; if $\dim M =2$ and $M$ is orientable, \cite{gafa} rules out the unit ball being a polyhedron. The situation is different when $\dim M =2$ and $M$ is not orientable :

\bigskip
\noindent {\bf Theorem C}
{\it Assume $M$ is a closed non-orientable surface. Then in every conformal class there exists a metric whose stable norm has a polyhedron as its unit ball.
}
\bigskip

Now we briefly describe the contents of the paper. In Sections 2 contains  basic facts  about non-orientable surfaces. In Sections 3 and 4 we have gathered prerequisites about  minimizing measures (in the sense of \cite{Mather91}) and stable norms. Some material from \cite{these} and \cite{gafa} has been included, either because it was not published, or because we found the redaction to be wanting. In Section 5 we prove the technical lemmas we need for our main theorems. Some consequences are derived, among which Proposition \ref{rational}, which says a minimizing measure with a rational homology class is supported on periodic orbits, and Lemma \ref{asymptote_fermee}, which says a geodesic asymptotic to a closed geodesic is not in the support of any minimizing measure.
In the last section we prove our main theorems.

\section{Preliminaries : non-orientable surfaces}

Let $(M,g)$ be a smooth, closed, non-orientable Riemannian manifold of dimension two.

\subsection{First homology group}

By the classical Surface Classification Theorem, any orientable closed surface is a connected sum of tori, any non-orientable closed surface is a connected sum of tori and projective planes. Since the connected sum of three projective planes is homeomorphic to the connected sum of a torus and a  projective plane, in fact any non-orientable surface is a connected sum of tori and one or two projective planes. Recall  that the connected sum of two projective planes is the Klein bottle $\K$.

Denote by $\Sigma_k$ an orientable surface of genus $k$ (that is, $\Sigma_k \cong \sharp^k \T^2$). We have
    \[ H_1( \Sigma_k \sharp \K, \R) \cong  H_1( \Sigma_k , \R) \oplus H_1(   \K, \R)\cong \R^{2k}\oplus \R
\]
whence the first Betti number $b_1(\Sigma_k \sharp \K)$ of $\Sigma_k \sharp \K$ is $2k+1$. Likewise,
\[
H_1( \Sigma_k \sharp \R P^2, \R) \cong  H_1( \Sigma_k , \R)\cong \R^{2k}
\]
and $
b_1(\Sigma_k \sharp \R P^{2})=2k$.

Similarly, we have
    \[ H_1( \Sigma_k \sharp \R P^{2}, \Z) \cong  H_1( \Sigma_k , \Z) \oplus H_1(  \R P^{2} , \Z)\cong \Z^{2k} \oplus \Z/2\Z
\]
and
    \[ H_1( \Sigma_k \sharp \K, \Z) \cong  H_1( \Sigma_k , \Z) \oplus H_1(   \K, \Z)\cong \Z^{2k}\oplus \Z \oplus \Z/2\Z.
\]
For any manifold $M$, the torsion-free part of $H_1( M, \Z)$ embeds as a lattice $\Lambda$ in $H_1( M, \R)$. We say
\begin{itemize}
    \item an element of $H_1( M, \R)$ is {\it integer} if it belongs to $\Lambda$
    \item a subspace of $H_1( M, \R)$ is {\it integer} if it is generated by integer classes
    \item  an element $h$ of $H_1( M, \R)$ is {\it rational} if $rh$ belongs to $\Lambda$ for some real number $r$.
\end{itemize}


\subsection{Orientation cover}

Let $\pi \co M_o \longrightarrow M$ be the orientation cover of
$M$. Then $M_o$ is an orientable surface endowed with a
fixed-point free, orientation-reversing involution $I$. Let
$I_{\ast}$ be the  involution of $H_1 (M_o,\R)$ induced by $I$,
and let $E_1$  (resp $E_{-1}$) be the eigenspace of $I_{\ast}$ for
the eigenvalue $1$ (resp$-1$). First observe  that

\begin{proposition}
 $E_1$  and $E_{-1}$ are Lagrangian for the (symplectic) intersection form $\inter$ on $H_1
(M_o,\R)$.
\end{proposition}

\proof Take $x,y \in E_1$  (resp $E_{-1}$). We have
$\inter (I_{\ast}(x), I_{\ast}(y) ) = \inter (x,y)$ but on the
other hand, since $I$ reverses the orientation of $M_o$, $\inter
(I_{\ast}(x), I_{\ast}(y) ) = -\inter (x,y)$ whence $\inter
(x, y) =0$, which proves that $E_1$  (resp
$E_{-1}$) is isotropic.
In particular $\dim E_1 \leq 2\inv b_1(M_o)$ and $\dim E_{-1} \leq 2\inv b_1(M_o)$.  Now since
$I_{\ast}$ is a linear involution, $\dim E_1 + \dim E_{-1}=b_1(M_o) $ whence $\dim E_1= \dim E_{-1}=2\inv b_1(M_o)$
that is, $E_1$  (resp $E_{-1}$) is Lagrangian for the symplectic
form $\inter$. \qed

\medskip

Furthermore

\begin{proposition}
$\ker \pi_{\ast}=E_{-1}$.
\end{proposition}

\proof Let $\gamma$ be a 1-cycle
in $M_o$ such that $ \pi_{\ast} ([\gamma])=0$. That is, $\pi
(\gamma)$ bounds a 2-chain $C$ in $M$. Then $\pi\inv (\pi
(\gamma))$ bounds the 2-chain $\pi\inv (C)$ in $M_o$. But $\pi\inv
(\pi (\gamma))=\gamma \cup I(\gamma)$, so $[\gamma]
+[I(\gamma)]=0$.
Conversely, if $\gamma$ is a 1-cycle
in $M_o$ such that $[\gamma]
+[I(\gamma)]=0$, then $\gamma$ and $I(\gamma)$ together bound a two-chain $C$ in $M_o$, so $\pi (\gamma)=\pi (I(\gamma))$ bounds the two-chain $\pi(C)$ in $M$, thus  $\left[\pi (\gamma)\right]=0$ in $H_1 (M,\R)$.
 \qed

\medskip

Consequently $\pi_{\ast}$
identifies $H_1 (M,\R)$ with $E_1$.

\section{Preliminaries : minimizing measures and stable norm}
The material of this section is taken from \cite{these} and was not published. Most of the ideas therein were presented to the second author by Albert Fathi.

Because we like our problems with some compacity, we introduce an alternative definition of the stable norm. It relies on invariant measures of the geodesic flow and is inspired by Mather's theory for Lagrangian systems. Then minimizing objects, in the form of measures (or asymptotic cycles as in \cite{Schwartzman}), exist
in every homology class. The question of whether a minimizing cycle exists becomes "are minimizing measures supported on closed geodesics ?".

Let $(M,g)$ be a compact Riemannian manifold  of any dimension with first Betti number $b_1(M) >0$.
Denote by
\begin{itemize}
  \item
  $T^{1}M $ the unit tangent bundle of $(M,g)$
    \item
    $p$ the canonical projection $T^{1}M \longrightarrow M$
    \item
    $\phi_t$ the geodesic flow in $T^{1}M$
\end{itemize}

\subsection{Minimizing measures.} Define $\mathcal{M}$ as the set of all probability measures on $T^{1}M$, endowed with the weak$^{\ast}$ topology.
 Then $\mathcal{M}$ is compact and metrizable (\cite{Dieudonne}, 13.4.2). Besides it embeds homeomorphically as a convex subset of  the dual to the vector space $C^{0}(T^{1}M)$ of continuous functions on $T^1M$.
 Let $\mathcal{M}_g$ be the subset of $\mathcal{M}$ that consists of $\phi_t$-invariant measures. Then $\mathcal{M}_g$ is closed in $\mathcal{M}$, hence compact, and convex.
Fix  an element $\mu$ of $\mathcal{M}_g$. By \cite{Mather91}, for any $C^{1}$ function $f$  on $M$, we have
    \[\int_{T^{1}M}df(x).v\; d\mu (x,v)=0.
\]
Thus, if $\omega$ is a smooth closed one-form on $M$, the integral
\[\int_{T^{1}M}\omega_x(v) d\mu (x,v)
\]
only depends on the cohomology class of $\omega$. By duality this endows $\mu$ with a homology class : $\left[\mu\right]$ is the unique element of $H_1(M,\R)$ such that
    \[<\left[\mu\right],\left[\omega\right]>=\int_{T^{1}M}\omega d\mu
\]
for any smooth closed one-form $\omega$ on $M$.
The map
    \[
    \begin{array}{rcl}
    \left[.\right] \co \mathcal{M}_g & \longrightarrow & H_1(M,\R) \\
    \mu & \longmapsto & \left[\mu \right]
    \end{array}
\]
is continuous and affine, so the image $\mathcal{B}_1$ of $\mathcal{M}_g$ in $H_1(M,\R)$ is compact and convex.

\begin{proposition}
$\mathcal{B}_1$ is the unit ball of the stable norm.
\end{proposition}

\proof We first prove that $\mathcal{B}_1$ is the unit ball of some norm $N$.

Denote by
    \[\begin{array}{rcl}
    \mathcal{I} \co T^{1}M & \longrightarrow T^{1}M \\
    (x,v) & \longmapsto (x,-v)
    \end{array}
\]
the canonical involution of $T^{1}M$.
We have, for any $(x,v)$ in $T^1 M$,
    \[ \phi_t (x,-v)= \phi_{-t}(x,v)
\]
so if  $\mu$ is  in  $\mathcal{M}_g$, then $\mathcal{I}_{\ast}\mu$ is again in $\mathcal{M}_g$. Let $\omega$ be a smooth closed one-form on $M$. We have
    \begin{eqnarray*}
    <\left[\mathcal{I}_{\ast}\mu\right],\left[\omega\right]>
    & = & \int_{T^{1}M}\omega_x(v) d\mathcal{I}_{\ast}\mu (x,v)\\
    & = & \int_{T^{1}M}\omega_x(-v) d\mu (x,v)\\
    & = & -\int_{T^{1}M}\omega_x(v) d\mu (x,v)\\
    & = &  -<\left[\mu\right],\left[\omega\right]>
\end{eqnarray*}
whence
\[\left[\mathcal{I}_{\ast}\mu\right]=-\left[\mu\right],
\]
 so $\mathcal{B}_1$ is centrally symetric.

Now let us show that $\mathcal{B}_1$ contains the origin in its interior. Fix a basis $h_1,\ldots h_n$ of $H_1(M,\R)$ such that $h_1,\ldots h_n$ are integer elements  of $H_1(M,\R)$. Let $\gamma_1,\ldots \gamma_n$ be closed geodesics parametrized by arc length such that $\left[\gamma_i\right]=h_i,\ i=1,\ldots n$ and let $\mu_1,\ldots \mu_n$ be the probability measures defined by
    \[ \int_{T^{1}M}f(x,v) d\mu_i (x,v):=
    \frac{1}{l_g(\gamma_i)}\int^{l_g(\gamma_i)}_{0}f\left(\gamma_i (t),\dot{\gamma_i}(t)\right)dt,\ i=1,\ldots n
\]
where $f \in C^0(T^1M)$.
We have, for any   smooth closed one-form $\omega$ on $M$
    \[<\left[\mu_i\right],\left[\omega\right]>=\frac{1}{l_g(\gamma_i)}<\left[\gamma_i\right],\left[\omega\right]>,\ i=1,\ldots n
\]

 whence
\[\left[\mu_i \right]=\frac{1}{l_g(\gamma_i)}\left[\gamma_i\right],\ i=1,\ldots n.
\]
Therefore $\mathcal{B}_1$ contains the points $\pm l_g(\gamma_i)\inv \left[\gamma_i\right],\ i=1,\ldots n$, so it contains their convex hull, which contains the origin in its interior because $\left[\gamma_i\right],\ i=1,\ldots n $ generate $H_1(M,\R)$.

So  $\mathcal{B}_1$ is the unit ball for some norm $N$ in $H_1(M,\R)$, which justifies the notation $\mathcal{B}_1$.

\medskip

Let us show this norm $N$ is no other than  the stable norm.
First we show $\|\ \| \geq N$.
Take $h$ in $H_1(M,\R)$ and $\epsilon >0$. Let $\sum_{i}r_i \gamma_i$ be a cycle such that the $\gamma_i$ are closed geodesics,
$\left[\sum_{i}r_i \gamma_i \right]=h$, and $\sum_{i}|r_i| l_g(\gamma_i)\leq \|h\| +\epsilon $. Reorienting the $\gamma_i$ if need be, we may assume that the $r_i$ are non-negative. Then the formula
    \[\int_{T^{1}M}f(x,v) d\mu (x,v):=
    \frac{\sum_i r_i \int^{l_g(\gamma_i)}_{0}f\left(\gamma_i (t),\dot{\gamma_i}(t)\right)dt}{\sum_{i}r_i l_g(\gamma_i)}
\]
defines an element of $\mathcal{M}_g$, with homology
    \[\left[\mu \right] = \frac{\left[\sum_{i}r_i \gamma_i \right]}{\sum_{i}r_i l_g(\gamma_i)}
\]
 By definition we have $N(\left[\mu \right])\leq 1$, whence, since $N$ is a norm
    \[ N\left( \left[\sum_{i}r_i \gamma_i \right] \right) \leq \sum_{i}|r_i| l_g(\gamma_i) \leq \|h\| +\epsilon.
\]
Since $\epsilon$ is arbitrary, we conclude that $\|\ \| \geq N$.

Now let us show that $\|\ \| \leq N$. It suffices to show that for any $\mu \in \mathcal{M}_g$, we have $\| [\mu]\| \leq 1$. Here we use the dual stable norm (see \cite{GLP}, 4.35). A norm $\|\ \|_0$ is defined on the space of $C^{1}$ closed one-forms on $M$ by
    \[ \|\omega \|_0 := \max \left\{ \omega_x(v) \co (x,v) \in T^1 M \right\}.
\]
This norm induces a norm on $H^{1}(M,\R)$ :  $\forall c \in H^{1}(M,\R)$,
    \[\|c \|_0 := \inf \left\{ \|\omega \|_0 \co \left[\omega\right]=c \right\}.
\]
\begin{lemma}[\cite{GLP}]
The norm $\|\  \|_0$ on $H^{1}(M,\R)$ is dual to the stable norm, that is, for any $h \in H_{1}(M,\R)$,
    \[ \|h  \| = \max \left\{ <c,h> \co c \in H^{1}(M,\R),\; \|c \|_0 \leq 1 \right\}.
\]
\end{lemma}
 In view of the above Lemma, what we need to show is that for any $c \in H^{1}(M,\R)$ such that $\|c \|_0 \leq 1$, we have  $<c,\left[\mu\right]> \leq 1$. As $\|c \|_0 \leq 1$, for all $\epsilon>0$ there exists a closed one-form $\omega$ such that $[\omega]=c$ and $|\omega_x (v)|\leq 1+\epsilon$ for all $(x,v) \in T^{1}M$. By the Ergodic Decomposition Theorem (\cite{Mane}, Theorem 6.4 p. 170) we have
    \[\int_{T^{1}M} \omega d\mu =  \int_{T^{1}M}\left\{\int_{T^{1}M} \omega d\mu_{x,v} \right\} d\mu (x,v)
\]
where, for $\mu$-almost every $(x,v)$,
    \[\int_{T^{1}M} \omega d\mu_{x,v}= \lim_{T \rightarrow +\infty} \frac{1}{T}\int^{T}_{0}\omega (\phi_t (x,v))dt.
\]
Since $\phi_t (x,v)$ is in $T^{1}M$ for all $t$, the above expression is $\leq 1+\epsilon$, which proves that $<\left[\omega\right],\left[\mu\right]> \leq 1+\epsilon$. Thus $<c,\left[\mu\right]> \leq 1$ for any $c \in H^{1}(M,\R)$ such that $\|c \|_0 \leq 1$ so $\|\ \| \leq N$.

Finally $\|\ \| = N$ and $\mathcal{B}_1=\{h \in H_1(M,\R) \co \|h\| \leq 1\}$. \qed

\bigskip

We say an element $\mu$ of $\mathcal{M}_g$ is {\it minimizing} if its homology class lies on the boundary $\mathcal{S}_1$  of $\mathcal{B}_1$, that is, if there exists a cohomology class $c$ such that $<c,\left[\mu\right]> =1$ and  $<c,h> \leq 1$ for all $h \in \mathcal{B}_1$.
\subsection{Link with Mather's theory.}  In this  paragraph we prove that the minimizing measures just defined are minimizing in the sense of Mather (\cite{Mather91}), which allows us to use Mather's Graph Theorem.

Consider the set $\mathcal{M}^{'}_{g}$  of all compactly supported, $\phi_t$-invariant probability measures on the tangent bundle $TM$ of $M$ and not just  $T^{1}M$ (here $\phi_t$ denote the geodesic flow in $TM$). We can define the homology class of an element of $\mathcal{M}^{'}_{g}$ just like we do for an element of $\mathcal{M}_g$. Mather's $\beta$-function is defined in \cite{Mather91} as
    \[  \begin{array}{rcl}
\beta \co H_1 (M,\R) & \longrightarrow & \R \\
h & \longmapsto &
\min \left\{\int_{TM}   \frac{1}{2}\|\cdot \|^{2}_{g}d\mu \co \mu \in \mathcal{M}^{'}_{g},\ \left[ \mu \right]=h \right\}
\end{array}
\]
where $\|(x,v)\|^{2}_{g}:=g_x(v,v)$ for all $(x,v) \in TM$.

The measures achieving the minimum for some $h$ are called {\it $h$-minimizing}. Next we show that this definition of minimizing agrees with ours.

\begin{proposition}
A minimizing measure $\mu \in \mathcal{M}_{g}$ is $[\mu]$-minimizing. Conversely an $h$-minimizing measure in $\mathcal{M}^{'}_{g}$ with $h \in \mathcal{S}_1$
is  in $\mathcal{M}_g$ ; in particular, it is minimizing. Furthermore, $2\beta = \|\ \|^{2}$.
\end{proposition}

\proof Let us begin by showing that $\beta$ is quadratic (i.e. $2$-homogeneous). Take $h \in H_1 (M,\R)$, $\mu$ an $h$-minimizing measure, $\lambda$ a real number. The formula
    \[f \longmapsto \int_{TM} f(x,\lambda v)d\mu(x,v)
\]
defines a probability measure on $TM$, whose homology class is $\lambda h$. Therefore we have
$ \beta (\lambda h) \leq \lambda^{2}\beta ( h)$ and likewise,
\[
 \beta ( h)= \beta (\frac{1}{\lambda} \lambda h)\leq \frac{1}{\lambda^{2}}\beta (\lambda h)
 \]
 whence $\beta (\lambda h) = \lambda^{2}\beta ( h)$.

 Now, since $2\beta$ and $ \|\ \|^{2}$ are both quadratic, proving that
    \[\mathcal{B}_1 = \left\{ h \in  H_1 (M,\R) \co \beta(h) \leq \frac{1}{2}\right\}
\]
suffices to prove that $2\beta = \|\ \|^{2}$. Note that
    \[\mathcal{B}_1 \subset \left\{ h \in  H_1 (M,\R) \co \beta(h) \leq \frac{1}{2}\right\}
\]
for if $\mu$ is an element of $\mathcal{M}_g$, we can view it as a measure on $TM$ supported on $T^{1}M$, thus
    \[\int_{TM}     \frac{1}{2}\|\cdot\|^{2}_{g}d\mu = \int_{TM}    \frac{1}{2}d\mu = \frac{1}{2}
\]
whence $\beta (\left[\mu \right])\leq 1/2$.

Conversely, let $h \in H_1 (M,\R)$ be such that $\beta (h)= 1/2$, and let $\mu$ be an $h$-minimizing measure. Then by \cite{Carneiro}
  the support of $\mu$ is contained in the energy level one half, that is, $T^{1}M$. Thus $\left[\mu \right] \in \mathcal{B}_1$, and
\[\mathcal{B}_1 = \left\{ h \in  H_1 (M,\R) \co \beta(h) \leq \frac{1}{2}\right\}\mbox{ whence } \beta = \frac{1}{2}\|\ \|^{2}.
\]

Besides, since $\beta (h)= 1/2$, $h$ lies on the boundary of $\mathcal{B}_1$, hence $\mu$ is minimizing in our sense.
Now we would like to prove that a minimizing measure in our sense  minimizes in the sense of Mather. Let $\mu \in \mathcal{M}_g$ be such that $\left[\mu \right] \in\mathcal{S}_1$. Then as we have just seen $\beta(\left[\mu \right])=1/2$ so
    \[\int_{TM}     \frac{1}{2}\|\cdot\|^{2}_{g}d\mu \geq \frac{1}{2}.
\]
On the other hand since $\mu$ is supported in $T^{1}M$, we have
\[\int_{TM}     \frac{1}{2}\|\cdot\|^{2}_{g}d\mu \leq \frac{1}{2}
\]
whence
\[\int_{TM}     \frac{1}{2}\|\cdot\|^{2}_{g}d\mu = \frac{1}{2}= \beta(\left[\mu \right])
\]
that is, $\mu$ is minimizing in the sense of Mather. \qed

\bigskip

The main reward of our efforts is that we may use Mather's Graph Theorem.
Let $c$ be a cohomology class such that    $<c,h> \leq 1$ for all $h \in \mathcal{B}_1$ and $<c,h_0> =1$ for some $h_0 \in \mathcal{B}_1$ (thus $\|c\|_0=1$).
We say a measure $\mu \in \mathcal{M}_g$ is {\it $c$-minimizing} if $<c,\left[\mu \right]> =1$.
Let $\mathcal{M}_c \subset T^{1}M$ be the union of the supports of all $c$-minimizing measures.
\begin{theorem}[Mather]
The restriction of $p$ to $\mathcal{M}_c$ is injective, and its inverse is Lipschitz.
\end{theorem}
This means that  minimizing measures can be identified with measured geodesic laminations.

\section{Flats of the unit ball}\label{flats}
Let $(M,g)$ be a closed Riemannian manifold of any dimension. We call
\begin{itemize}
    \item supporting subspace to the unit ball of the stable norm, any affine subspace of $H_1(M,\R)$ that meets the unit sphere but not the open unit ball
    \item flat of the unit ball, the intersection of the unit sphere with a supporting subspace
    \item dimension of a flat, the dimension of the affine subspace it generates in $H_1(M,\R)$
    \item interior of a flat, its interior in the affine subspace it generates.

\end{itemize}
 As a trivial example, all points of the unit sphere are zero-dimensional flats. If $c$ is a cohomology class of dual stable norm one, that is, $<c,h> \leq 1$ for all $h \in \mathcal{B}_1$ and $<c,h_0> =1$ for some $h_0 \in \mathcal{B}_1$, then
    \[ \left\{ h \in H_1(M,\R) \co <c,h>=1 \right\}
\]
is a supporting hyperplane to the unit ball of the stable norm, and
\[ \left\{ h \in  \mathcal{B}_1\co <c,h>=1 \right\}
\]
is a flat, which may or may not be trivial. Note that by the Hahn-Banach Theorem, any supporting subspace is contained in a supporting hyperplane. So for any flat $F$, there exists $c \in H^1 (M,\R)$ such that
\[<c,h> \leq 1 \   \forall h \in \mathcal{B}_1 \mbox{ and } F \subset \left\{ h \in  \mathcal{B}_1\co <c,h>=1 \right\}.
\]
By \cite{Mane92}, if a minimizing measure $\mu$ is ergodic, then $\left[\mu\right]$ is an extremal point of $\mathcal{B}_1$, hence it cannot be in the interior of any non-trivial flat. In particular non connected minimizing cycles, when they exist,  are the simplest examples of non-trivial flats.
Recall Proposition 4 of \cite{gafa} :
\begin{lemma}\label{extension_plat}
Let $F_1$ and $F_2$ be two flats of the unit ball, both containing a point $h_0$ such that $h_0$ is an interior point of $F_1$. Then there exists a flat $F$ containing $F_1 \cup F_2$.
\end{lemma}
\proof
Let $c \in H^1 (M,\R)$ be such that
\[ <c,h> \leq 1 \   \forall h \in \mathcal{B}_1 \mbox{ and } F_2 \subset \left\{ h \in  \mathcal{B}_1\co <c,h>=1 \right\}.
\]
Then $c$ restricted to the convex set $F_1$ has  a maximum at the interior point $h_0$. Since $c$ is linear, this implies that $c$ is constant on $F_1$. Hence
\[ \left\{ h \in  \mathcal{B}_1\co <c,h>=1 \right\}
\]
is a flat containing $F_1 \cup F_2$.
\qed

\begin{lemma}
Let $F_1$ and $F_2$ be two flats of the unit ball, both containing a point $h_0$ in their interiors. Then there exists a flat $F$ containing $F_1 \cup F_2$ such that  $h_0$ is an interior point of $F$.
\end{lemma}
\proof
Let
\begin{itemize}
    \item $V_i$, $i=1,2$ be the underlying vector space of the affine space generated by $F_i$
    \item $V := V_1 + V_2$
    \item $A$ be the affine subspace $h_0+V$
    \item $F := A \cap \mathcal{B}_1$
    \item $c \in H^1 (M,\R)$ be given by the previous Lemma such that
\[ <c,h> \leq 1 \   \forall h \in \mathcal{B}_1 \mbox{ and } F_1, F_2 \subset \left\{ h \in  \mathcal{B}_1\co <c,h>=1 \right\}.
\]
\end{itemize}
Since $F_1, F_2 \subset \left\{ h \in  \mathcal{B}_1\co <c,h>=1\right\}$ we have $A \subset \left\{ h \in  \mathcal{B}_1\co <c,h>=1\right\}$ so $A$ is a supporting subspace whence $F$ is a flat. Besides, since $F$ is convex and contains $F_1$ and $F_2$, it contains the convex hull $C$ of $F_1$ and $F_2$. Now, since $h_0$ is interior to both $F_1$ and $F_2$, there exist open neighborhoods of zero  $U_1$, $U_2$ in  $V_1$, $V_2$ respectively such that $h_0 + U_i \subset F_i$, $i=1,2$. So the convex hull of $h_0+U_1$ and  $h_0+U_2$ is open in $A$, and contained in $C$, hence in $F$. Thus $h_0$ is an interior point of $F$.
\qed

\medskip

The former Lemma means that for any homology class $h$, there exists an unique maximal flat containing $h$ in its interior.

\bigskip

\noindent{\bf Orientable surfaces.}

\medskip

Assume, for the remainder of this section, that $M$ is an orientable surface.  If $F$ is a flat of the unit ball of the stable norm,  $h_1, h_2 \in F$ and $\mu_1$, $\mu_2$ are minimizing measures such that $\left[\mu_i\right]=h_i, \  i= 1,2$, then by Mather's Graph Theorem the supports of $\mu_1$ and  $\mu_2$ do not intersect transversally so $\inter (h_1,h_2)=0$. Thus the vector space generated by $F$ in $H_1 (M,\R)$ is isotropic with respect to the symplectic intersection form. In particular its dimension is $\leq 2\inv b_1(M)$ so $\dim F \leq 2\inv b_1(M)-1$. This was first observed by M.J. Carneiro (\cite{Carneiro}). In \cite{gafa} this upper bound is proved to be optimal, so non-trivial flats always exists for orientable surfaces of genus $\geq 2$ (see Theorem \ref{local} and its corollary).

It is proved in \cite{gafa} (Proposition  6, which surprisingly we have seen no need to reprove) that a flat containing a rational point in its interior is a finite polyhedron with  at most $3(2\inv b_1(M)-1)$ vertices. Furthermore the vertices are rational homology classes which have connected minimizing cycles.

\section{Technical lemmas}

\subsection{Key lemma, one-sided case}
After writing this lemma we came across reference \cite{Scharlemann}, where a similar result is proved in a topological setting.
Lemma \ref{lemmanonor} and its orientable companion Lemma \ref{key_o} are improved versions  of Lemmas 15, 16, 17 of \cite{gafa}.
The purpose of Lemma \ref{lemmanonor} was to mimic the approach of \cite{gafa}. Later on, inspired by \cite{Fathi} we realized it is simpler to use the orientation cover. So Lemma \ref{lemmanonor} won't be used in the proofs of our main theorems. Still we believe it is interesting in its own right. We do use it, however, in the proof of Lemma \ref{asymptote_fermee}.

Let $\gamma_1$ be a closed, one-sided, simple geodesic on a non-orientable surface $M$.
\begin{lemma}\label{lemmanonor}
There exists a neighborhood $V_1$ of $(\gamma_1,\dot\gamma_1)$ in $T^{1}M$ such that, for  any simple  geodesic $\gamma$, if $(\gamma,\dot\gamma)$ enters (resp. leaves) $V_1$  then $\gamma$ is forever trapped in $p(V_1)$ in the future (resp. past), that is
    \[\exists t_1 \in \R,\; \forall t \geq (\mbox{resp }\leq ) t_1, \; \gamma(t)\in p(V_1).
\]
\end{lemma}
\proof Let $U_1$ be a neighborhood of $\gamma_1$ in $M$ homeomorphic to a M\H{o}bius strip.
Let $P:=\gamma_1 (0)$ be a point of $\gamma_1$, and let $\delta$ be an smooth open arc transverse at $P$ to $\gamma_1$, such that $U_1 \setminus \delta$ is simply connected. Let $V_1$ be the neighborhood of $(\gamma_1,\dot{\gamma_1})$ in  $T^{1}M$ defined by
\begin{enumerate}
    \item
    $p(V_1)=U_1$
    \item
    $\forall (x,v) \in V_1, p(\phi_t (x,v))$ intersects transversally $\delta$ at least  three times $t_1 <t_2<t_3$ (the points
    $p(\phi_{t_i} (x,v)), i=1,2,3$ may coincide if $p(\phi_t (x,v))$ is a closed geodesic)
    \item
    we would like a condition along the lines of "the  geodesics that enter $V$ always cross $\delta$ in the same direction as $\gamma_1$". This takes some precaution because $M$ is not orientable. So we choose a smooth vector field $X$ in $U_1$, transverse to $\delta$, which has $\gamma_1$ as a trajectory and such that every other trajectory is closed and homotopic to the boundary of $U_1$, that is, bounds a  M\H{o}bius strip containing $\gamma_1$. We require that
    \[\forall (x,v) \in V_1, g(X(x),v) >0.
\]

\end{enumerate}
These conditions define an open set of $T^{1}M$ because $\delta$ is an open arc and we demand that the intersections be transverse.

Now consider a simple geodesic $\gamma$ such that for some $t \in \R$, $(\gamma(t),\dot\gamma(t))\in V_1$. Let $t_1 < t_2 <t_3 $ be such that $\gamma(t_1),\gamma(t_2),\gamma(t_3)$ are the transverse intersection points of $\gamma$ and $\delta$ given by the definition of $V_1$.
\subsubsection{First case}
\begin{figure}[h]
\leavevmode \SetLabels
\L(.46*.4) $\gamma_1$\\
\L(.46*.85) $\gamma$\\
\L(.2*.55) $\delta$\\
\L(.16*.8) $\gamma(t_1)$\\
\L(.16*.9) $\gamma(t_3)$\\
\L(.16*.28) $\gamma(t_2)$\\
\L(.78*.7) $\gamma(t_2)$\\
\L(.78*.2) $\gamma(t_1)$\\
\L(.78*.1) $\gamma(t_3)$\\
\L(.62*.6) $U'_1$\\
\endSetLabels
\begin{center}
\AffixLabels{\centerline{\epsfig{file =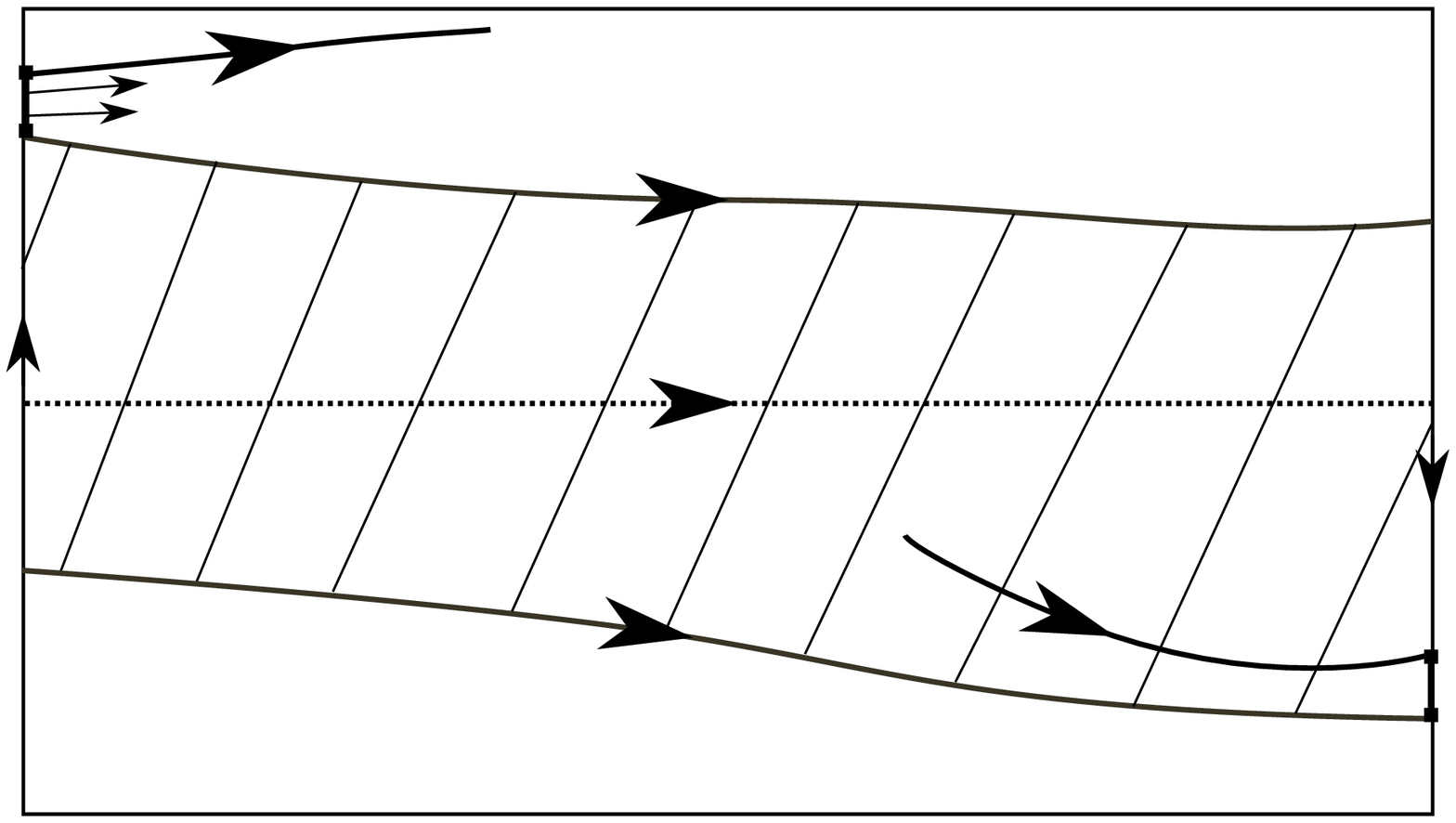,width=7cm,angle=0}}}
\end{center}
\caption{}
\label{casnonor}
\end{figure}

Assume $\gamma(t_3)$ is farther away from $P$ than  $\gamma(t_1)$ with respect to the distance on $\delta$ induced by the metric of $M$. The domain  $U'_1$ bounded by $\gamma(\left[ t_1,t_3 \right])$  and the subsegment of  $\delta$ joining $\gamma(t_1)$ with $\gamma(t_3)$ is homeomorphic to a M\H{o}bius strip. The geodesic $\gamma$ does not self-intersect, hence it can only cut the boundary of  $U'_1$ along $\delta$. By Condition 3 of the definition of $V_1$, $\gamma$
can only intersect $\delta$ from left to right as pictured in Figure \ref{casnonor}, that is, outwards of $U'_1$. Therefore $\gamma$ is trapped in $U'_1$ in the past.

\subsubsection{Second case}
Assume $\gamma(t_1)$ is farther away from $P$ than  $\gamma(t_3)$ with respect to the distance on $\delta$ induced by the metric of $M$. We prove as in the first case that $\gamma$ is trapped in $U'_1$ in the future.

\subsubsection{Third case}
Assume $\gamma(t_1) = \gamma(t_3)$. Then, since $\gamma$ doesn't self-intersect, it must be a closed geodesic, hence is trapped in $U_1$  in both past and future.
\qed

\subsection{Key lemma, two-sided case}
Let $\gamma_2$ be a closed, two-sided, simple geodesic on a surface $M$, orientable or not. Again, we shall only use the orientable case here, but the extra generality comes for free. Let $U_2$ be a neighborhood of $\gamma_2$ in $M$ homeomorphic to an annulus. Choose a symplectic form $\omega$ in $U_2$, yielding a local orientation of $U_2$.
\begin{lemma}\label{key_o}
 There exists a neighborhood $V_2$ of $(\gamma_2,\dot\gamma_2)$ in $T^{1}M$ such that any simple  geodesic $\gamma$, if  $(\gamma,\dot\gamma)$ enters (resp. leaves) $V_2$  then either $\gamma$ intersects $\gamma_2$ or $\gamma$ is forever trapped in $p(V_2)$ in the future (resp. past), that is
    \[\exists t_2 \in \R,\; \forall t \geq (\mbox{resp }\leq ) t_2, \; \gamma(t)\in p(V_2).
\]
Besides, all intersections with $\gamma_2$ have the same sign with respect to $\omega$.
\end{lemma}
\proof
Let $P:=\gamma_2 (0)$ be a point of $\gamma_2$, and let $\delta$ be an smooth open arc transverse at $P$ to $\gamma_2$, such that $U_2 \setminus \delta$ is simply connected.
 Assume $\delta$ is oriented so that
$\omega(P)(\dot\gamma_2,\dot\delta)>0$. Let $V_2$ be the neighborhood of $(\gamma_2,\dot{\gamma_2})$ in  $T^{1}M$ defined by
\begin{enumerate}
    \item
    $p(V_2)=U_2$
    \item
    $\forall (x,v) \in V_2,\  p(\phi_t (x,v))$ intersects transversally $\delta$ at least twice before intersecting $\gamma_2$, if it  intersects $\gamma_2$ at all ; and if it does, it must intersect transversally $\delta$ at least twice more before either leaving $U_2$ or meeting $\gamma_2$ again
    \item
    the  geodesics that enter $V_2$ always cross $\delta$ in the same direction as $\gamma_2$, that is,
    \[
    \forall x \in \delta,\; \forall v \in T^{1}_{x}M \mbox{ such that } (x,v) \in V_2,\; \omega(x)(v,\dot\delta (x))>0.
    \]
\end{enumerate}
These conditions define an open set of $T^{1}M$ because $\delta$ is an open arc and we demand that the intersections be transverse.

Now consider a simple geodesic $\gamma$ such that for some $t \in \R$, $(\gamma(t),\dot\gamma(t))\in V_2$. Let $t_1 < t_2  $ be such that $\gamma(t_1),\gamma(t_2)$ are the first two transverse intersection points of $\gamma$ and $\delta$ given by the definition of $V_2$.

\subsubsection{First case}
\begin{figure}[h]
\leavevmode \SetLabels
\L(.47*.4) $\gamma_2$\\
\L(.45*.82) $\gamma$\\
\L(.18*.3) $\delta$\\
\L(.18*.5) $P$\\
\L(.15*.7) $\gamma(t_1)$\\
\L(.15*.84) $\gamma(t_2)$\\
\L(.8*.84) $\gamma(t_2)$\\
\L(.8*.7) $\gamma(t_1)$\\
\L(.3*.6) $U'_2$\\
\endSetLabels
\begin{center}
\AffixLabels{\centerline{\epsfig{file =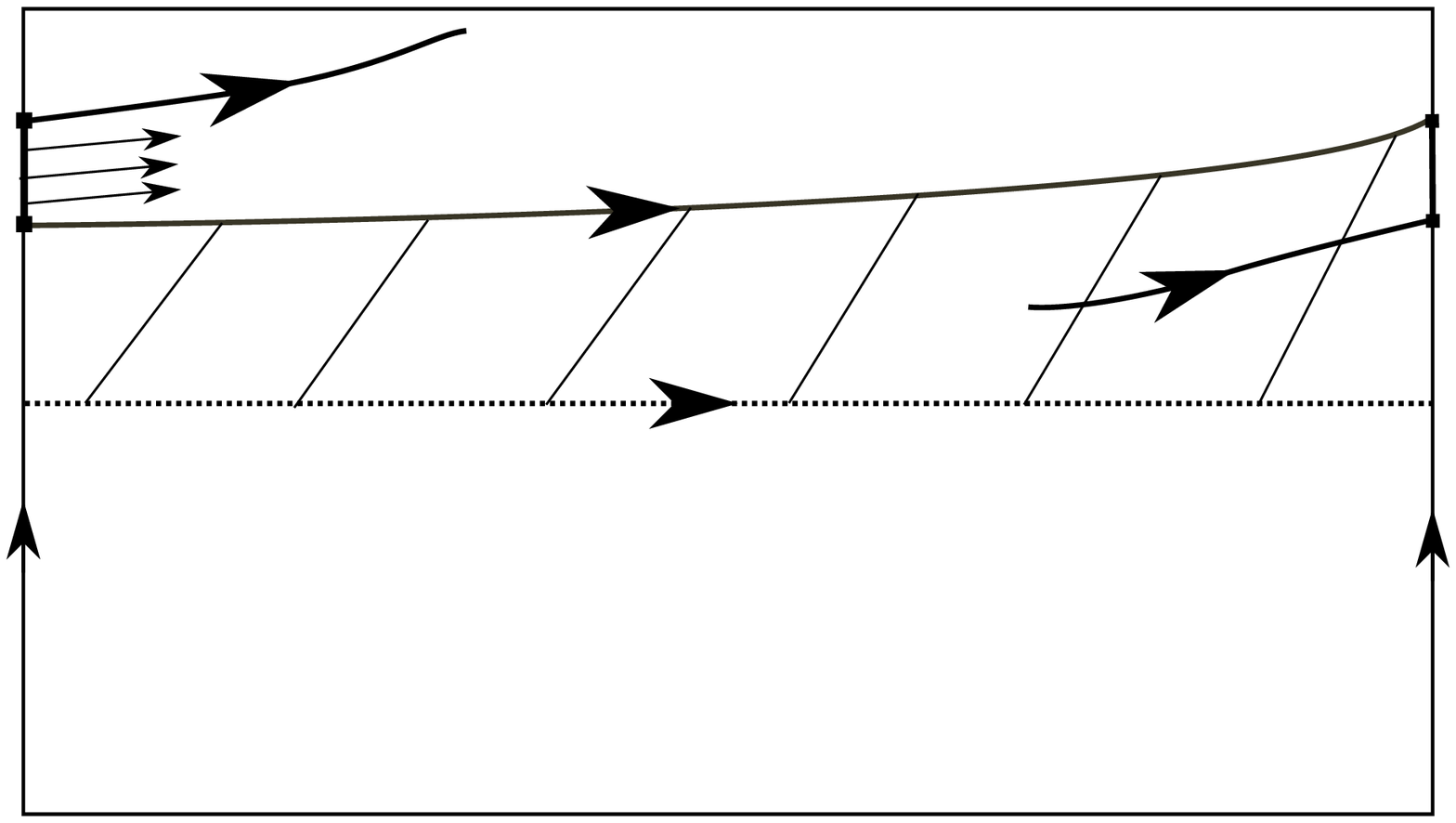,width=7cm,angle=0}}}
\end{center}
\caption{}
\label{casor}
\end{figure}

Assume $\gamma(t_2)$ is farther away from $P$ than  $\gamma(t_1)$ with respect to the distance on $\delta$ induced by the metric of $M$. The domain  $U'_2$ bounded by $\gamma(\left[ t_1,t_2 \right])$  and the subarc of  $\delta$ joining $\gamma(t_1)$ with $\gamma(t_2)$ on one side, and by $\gamma_2$ on the other side is homeomorphic to an annulus. The geodesic $\gamma$ is simple so it cannot self-intersect, hence it can only cut the boundary of  $U'_2$ along $\delta$ or $\gamma_2$. By Condition 3 of the definition of $V_2$, $\gamma$
can only intersect $\delta$ from left to right as pictured in Figure \ref{casor}, that is, outwards of $U'_2$. Therefore $\gamma$ either intersects $\gamma_2$ or is trapped in $U'_2$ in the past.

\subsubsection{Second case}
Assume $\gamma(t_1)$ is farther away from $P$ than  $\gamma(t_2)$ with respect to the distance on $\delta$ induced by the metric of $M$. Likewise we prove that $\gamma$ either intersects $\gamma_2$ or is trapped in $U'_2$ in the future.

\subsubsection{Third case}
Assume $\gamma(t_1) = \gamma(t_2)$. Then, since $\gamma$ doesn't self-intersect, it must be a closed geodesic, and the conclusion  readily follows.

\medskip

\begin{figure}[h]
\leavevmode \SetLabels
\L(.47*.4) $\gamma_2$\\
\L(.45*.82) $\gamma$\\
\L(.55*.58) $\gamma(0)$\\
\L(.18*.5) $\delta$\\
\L(.15*.2) $\gamma(t_4)$\\
\L(.15*.35) $\gamma(t_3)$\\
\L(.8*.84) $\gamma(t_1)$\\
\L(.8*.7) $\gamma(t_2)$\\
\L(.29*.57) $U''_2$\\
\endSetLabels
\begin{center}
\AffixLabels{\centerline{\epsfig{file =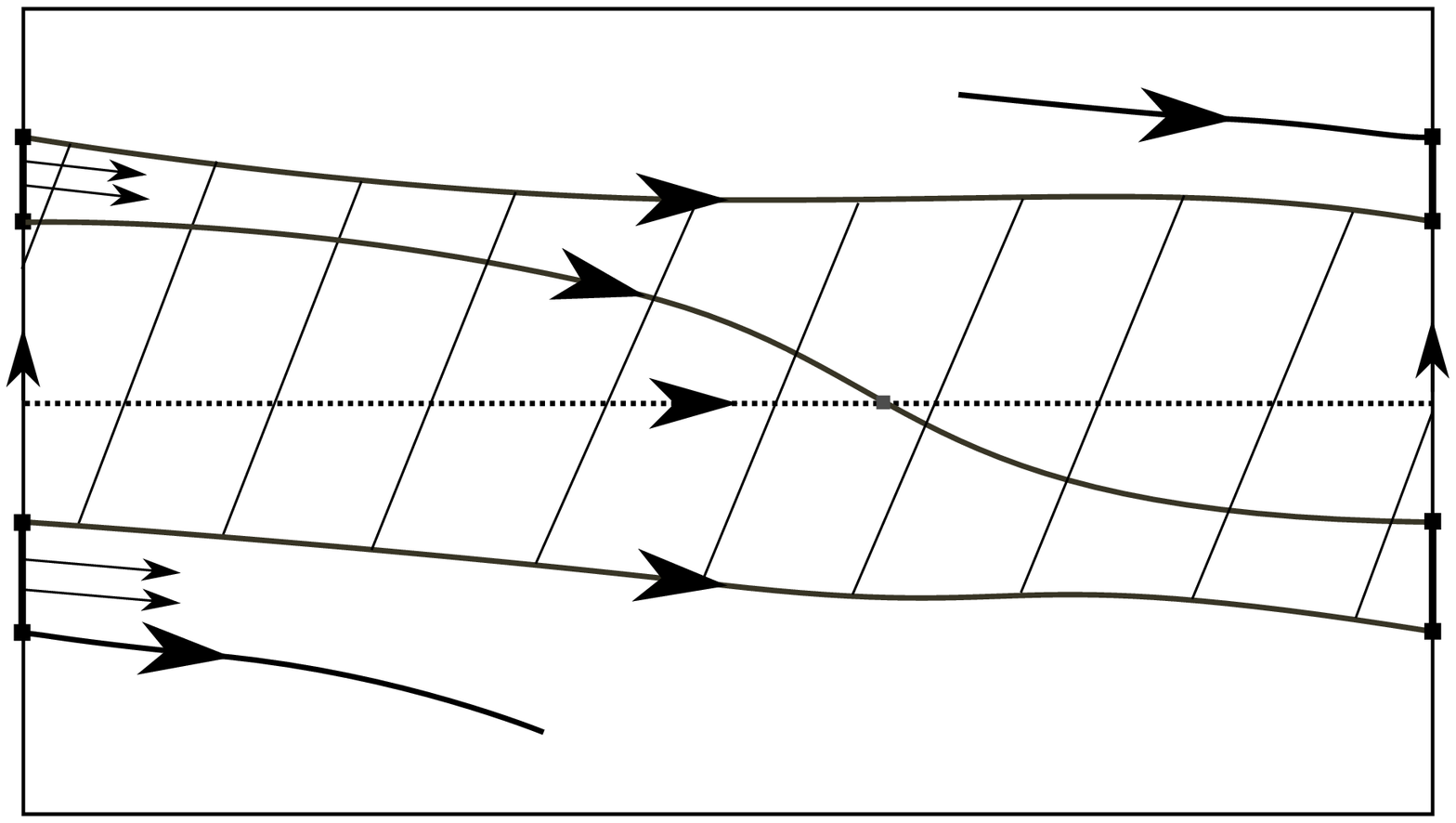,width=7cm,angle=0}}}
\end{center}
\caption{}
\label{casor3}
\end{figure}

We still have to prove the statement about the sign of the intersections. Assume $\gamma$ cuts $\gamma_2$ once with  positive sign, that is, downwards in Figure \ref{casor3}. Assume for convenience that the intersection point is $\gamma (0)$. Let
$t_1 < t_2 <0< t_3 <t_4 $ be such that $\gamma(t_1),\gamma(t_2)$ are the last two transverse intersection points of $\gamma$ and $\delta$ before $\gamma$ meets $\gamma_2$, and $\gamma(t_3),\gamma(t_4)$ are the first two transverse intersection points of $\gamma$ and $\delta$ after $\gamma$ meets $\gamma_2$.  The domain  $U''_2$ bounded by $\gamma(\left[ t_1,t_2 \right])$  and the subarc of  $\delta$ joining $\gamma(t_1)$ with $\gamma(t_2)$ on one side, and by $\gamma(\left[ t_3,t_4 \right])$  and the subarc of  $\delta$ joining $\gamma(t_3)$ with $\gamma(t_4)$ on the other side,  is homeomorphic to an annulus and contains $\gamma_2$ in its interior.

The geodesic $\gamma$ does not self intersect so it cannot enter $U''_2$ through segments of $\gamma$. It only intersects $\delta$ from left to right, that is, either between $\gamma(t_3)$ and $\gamma(t_4)$ and outwards of $U''_2$, or between $\gamma(t_1)$ and $\gamma(t_2)$ and inwards of $U''_2$. So it can only  enter $U''_2$ through $\delta$  between $\gamma(t_1)$ and $\gamma(t_2)$, that is, from above in Figure \ref{casor3}. Therefore it always cut  $\gamma_2$ with positive sign.
\qed

\subsection{Consequences of the key lemmas}

For the ease of the reader, rather than loading up the sentences with``resp."  we have split the next proposition in two, one part for the one-sided case and the other for the two-sided case.
\begin{proposition}\label{one-sided}
Let $\gamma_{1}$ be a closed, simple, one-sided geodesic on a surface $M$, orientable or not. There exists a neighborhood $V_1$ of $\left(\gamma_1,\dot\gamma_1\right)$ in $T^{1}M$ such that for any simple geodesic $\gamma$, if $(\gamma, \dot\gamma)$ enters (resp. leaves) $V_1$, then
\begin{itemize}
    \item
     $\gamma$ is a closed geodesic homotopic to $\gamma_1$ or $\gamma_1.\gamma_1$, the latter being meant as a product in $\pi_1 (M)$
    \item
    or $\gamma$ is positively (resp. negatively) asymptotic to a closed geodesic homotopic  to  $\gamma_1$ or $\gamma_1.\gamma_1$.
\end{itemize}

\end{proposition}

\begin{proposition}\label{two-sided}
Let $\gamma_2$ be a  closed, simple, two-sided geodesic. There exists a neighborhood $V_2$ of $\left(\gamma_2,\dot\gamma_2\right)$ in $T^{1}M$ such that for any simple geodesic $\gamma$, if $\gamma$ enters (resp. leaves) $V_2$, then
\begin{itemize}
  \item
  either  $\gamma$ is  a closed geodesic homotopic to $\gamma_2$
    \item
or  $\gamma$ is asymptotic to a closed geodesic homotopic to $\gamma_2$
    \item
    or $\gamma$ intersects $\gamma_2$,   and all intersections  have the same sign with respect to some  orientation of $p (V_2)$.
\end{itemize}
\end{proposition}
\proof
Let us  prove Proposition \ref{two-sided}. Let $V_1$ be a neighborhood  of $\left(\gamma_1,\dot\gamma_1\right)$ in $T^{1}M$ given by Lemma \ref{key_o}, and small enough so it does not contain any contractible closed geodesic.  Let $\gamma$ be a simple geodesic such that $(\gamma, \dot\gamma)$ enters  $V_1$. Let  $t_1$ be such that $\left(\gamma, \dot{\gamma}\right)\left(\left[t_1,+\infty \right[ \right)\subset V_1$. Extend $\dot{\gamma}\left(\left[t_1,+\infty \right[ \right)$ to a smooth vector field  in the annulus $p(V_1)$.
Since the annulus may be embedded in the two-sphere,   the Poincar\'e-Bendixon Theorem applies. So $\gamma$ is either  a fixed point, a cycle of fixed points and heteroclinic orbits, or a closed orbit, or asymptotic to one of the former. Four out of six   cases are impossible here because $\gamma$ is a geodesic so its velocity is constant, hence cannot go to zero. Besides, a closed orbit of a vector field must be a simple closed curve, and a non-contractible simple closed curve in an annulus is homotopic to  the boundary of the annulus.
 This proves Proposition \ref{two-sided}. The proof of Proposition \ref{one-sided} is identical, mutatis mutandis.
\qed

\subsection{Geodesics asymptotic to closed geodesics}
\begin{lemma}\label{asymptote_fermee}
If a  geodesic $\gamma$ is asymptotic to a simple closed geodesic, then $(\gamma,\dot{\gamma})$ is not in  the support of any minimizing measure.
\end{lemma}
\proof
Let
\begin{itemize}
    \item
    $\gamma_0$ be a simple closed geodesic
    \item
    $\gamma$ a geodesic asymptotic to $\gamma_0$
    \item
    $V_0$ be a neighborhood  of $\left(\gamma_0,\dot\gamma_0\right)$ in $T^{1}M$ given by Lemma \ref{lemmanonor} or \ref{key_o}
    depending on whether $\gamma_0$ is one-sided or two-sided, and such that $\left(\gamma(0),\dot\gamma(0)\right) \notin V_0$

    \item
    $V$ a neighborhood of $\left(\gamma(0),\dot\gamma(0)\right)$ such that $V $ is disjoint from $V_0$ but for large enough $t$, $\phi_t(V)\subset V_0$
    \item
    $\mu$ a minimizing measure.
\end{itemize}
Assume $(\gamma,\dot{\gamma})$ is contained in the support of $\mu$. Then since $\supp (\mu)$ is closed and invariant under the geodesic flow, it contains the $\alpha$-and $\omega$-limit sets of $(\gamma,\dot{\gamma})$, in particular it contains $\left(\gamma_0,\dot\gamma_0\right)$. Besides, $\mu(V)>0$. By the Ergodic Decomposition Theorem (\cite{Mane}, Theorem 6.4 p. 170), we have
    \[\mu(V)= \int_{T^{1}M}\mu_{x,v}(V)d\mu(x,v)
\]
where, for $\mu$-almost every $x,v$,
    \[\mu_{x,v}(V)= \lim_{T \rightarrow +\infty} \frac{1}{T}\int^{T}_{0}\chi_V(\phi_t(x,v))dt
\]
denoting by $\chi_V$ the characteristic function of $V$. Thus for some $(x,v)$ we have $\mu_{x,v}(V)>0$.
So for some $t$ in $\R$, $\phi_{t}(x,v) \in V$. By our hypothesis on $V$, this implies that for $t$ large enough $\phi_{t}(x,v) \in V_0$. But since $\left(\gamma_0,\dot\gamma_0\right)$ is contained in the support of $\mu$, by Mather's Graph Theorem $p(\phi_{t}(x,v))$ cannot intersect $\gamma_0$. Thus by Proposition \ref{one-sided} or  \ref{two-sided}  the geodesic $p(\phi_{t}(x,v))$ is asymptotic to a geodesic homotopic to $\gamma_0$ whose lift to $T^1 M$ is contained in $V_0$. Therefore $\phi_{t}(x,v)$ never comes back to $V$, whence
\[ \lim_{T \rightarrow +\infty} \frac{1}{T}\int^{T}_{0}\chi_V(\phi_t(x,v))dt=0
\]
which contradicts the fact that $\mu_{x,v}(V)>0$.
\qed
\subsection{Minimizing measures with rational homology classes}
The  Proposition below was stated as Lemma 2.1.6 in \cite{these} and Proposition 5 of \cite{gafa}, although the announcement of its proof was greatly exaggerated.

\begin{proposition}\label{rational}
Let $M$ be a closed surface, orientable or not, with a Riemannian metric.
If $h$ is a rational homology class and $\mu$ is an $h$-minimizing measure, then the support of $\mu$ consists of periodic orbits.
\end{proposition}

First we need a
\begin{lemma}\label{relevee_minimisante}
Let $M$ be a closed non-orientable surface with a Riemanninan metric $g$ and $\pi : (M_o,\tilde{g}) \rightarrow (M,g)$ denote the Riemannian orientation  cover.

If $\mu \in \mathcal{M}_{g}$ is a $c$-minimizing measure where $c$ is a cohomology class with $\|c\|_0=1$, then there exists $\nu \in \mathcal{M}_{\tilde{g}}$ such that $\pi_{\ast}(\nu)= \mu$, $\nu$ is  $I_{\ast}$-invariant and $\pi^{\ast}(c)$-minimizing.
\end{lemma}

\proof[Proof of the Lemma] Let us assume for the time being that $\mu$ is ergodic. That is, there exists $(x,v) \in T^{1}M$  such that  for any continuous function $F$ on  $T^{1}M$,
    \[\int_{T^{1}M}F d\mu =  \lim_{T\rightarrow +\infty} \frac{1}{T}\int^{T}_{0} F(\phi_t(x,v))dt
\]
Let us lift the orbit $\phi_t(x,v)$ to an orbit $\tilde{\phi}_{t}(x,v)$ of the geodesic flow of $(M_o,\tilde{g}) $.
Let $\nu_T$ be the probability measure on $T^{1}M_o$ defined by
    \[\int_{T^{1}M_o}F d\nu_T=   \frac{1}{T}\int^{T}_{0} F(\tilde{\phi}_t(x,v))dt
\]
for any continuous function $F$ on  $T^{1}M_o$.
Since the set of probability measures  on $T^{1}M_o$ is compact for the weak$^{\ast}$ topology, there exists a sequence $T_n \rightarrow +\infty$ such that $\nu_{T_n}$ converges to some measure $\nu$. Then $\nu$ is invariant by the geodesic flow on $M_o$, that is, $\nu \in \mathcal{M}_{\tilde{g}}$. Besides   $\pi_{\ast}(\nu)= \mu$ since, for any continuous function $F$ on  $T^{1}M$
\begin{eqnarray*}
\int_{T^{1}M}F d\pi_{\ast}(\nu) &=&     \int_{T^{1}M_o}F\circ \pi\;  d\nu =
\lim_{n \rightarrow +\infty} \frac{1}{T_n}\int^{T_n}_{0} F(\pi (\tilde{\phi}_t(x,v)))dt \\
&=&  \lim_{n \rightarrow +\infty} \frac{1}{T_n }\int^{T_n }_{0} F(\phi_t(x,v))dt = \int_{T^{1}M}F d\mu.
\end{eqnarray*}
Furthermore, since $\mu$ is a $c$-minimizing measure, for all $\epsilon >0$ there exists a closed one-form $\omega$ such that $\left[\omega\right]=c$,
$|\omega_x (v)| \leq 1+\epsilon$ for all $x \in M$, $v \in T^{1}_{x} M$, and $\int \omega d \mu =1$. Set $\tilde{\omega}:= \pi^{\ast} (\omega)$, then  $\tilde c=[\tilde\omega]=\pi^\ast(c)$ and $|\tilde\omega_x (v)| \leq 1+\epsilon$ for all $x \in M_o$, $v \in T^{1}_{x} M_o$, and $\int \tilde{\omega}d\nu = 1$. So $\nu$ is $\pi^\ast (c)$-minimizing.  Notice that $\tilde{\omega}$ is $I^{\ast}$-invariant, so $I_{\ast}\nu$ is also $\pi^\ast (c)$-minimizing. Then  so is $2^{-1}(\nu + I_{\ast}\nu)$, which is $I_{\ast}$-invariant. This proves the ergodic case of the Lemma as $\pi_\ast (I_\ast \nu)=\mu$.

Now consider the map $\pi_{\ast}$ between the two compact convex sets
\[
\left\{ \mu \in \mathcal{M}_{\tilde{g}} \co I_{\ast} (\mu) = \mu \right\}
\]
and  $\mathcal{M}_{g}$.
It is affine and surjective onto the extremal points of $\mathcal{M}_{g}$, hence surjective onto $\mathcal{M}_{g}$.
\qed

\bigskip
\proof[Proof of the Proposition]
\noindent First let us address the case when $M$ is orientable.

Let $h$ be a rational homology class and $\mu$ be an $h$-minimizing measure. Then $\inter (h, H_1(M,\Z))$ is a discrete subgroup of $\R$. Assume the projection $p(\supp \mu)$ of the support of $\mu$ to $M$ contains a non-closed geodesic $\gamma$. Since $M$ is compact $\gamma$ has a limit point, say $x_{\gamma}$ in $M$. Let $t_n$ be an increasing sequence of real numbers such that $\gamma(t_n) \longrightarrow x_{\gamma}$ when $n \longrightarrow \infty$. Denote by $\gamma_n$ the closed curve obtained by closing up $\gamma (\left[t_n,t_{n+1}\right])$ with a geodesic segment $\delta_n$ of length $d(\gamma(t_n),\gamma(t_{n+1}))$. Such a segment is unique for $n$ large enough because $d(\gamma(t_n),\gamma(t_{n+1}))$ tends to zero. We claim that
$\inter (h, \left[\gamma_n \right]))$ is not zero for $n$ large enough, and tends to zero, which contradicts the discreteness of $\inter (h, H_1(M,\Z))$.

By Mather's Graph Theorem (\cite{Mather91}), for any $x$ in $p (\supp \mu )$, there exists a unique geodesic, denoted $\gamma_{x}$, which is the projection of an orbit in $\supp \mu$ and such that $\gamma_{x}(0)=x$. To clear up the  notations  we denote by $\gamma_{\gamma}$
the orbit $\gamma_{x}$ with $x=x_{\gamma}$. Call
    \[ R_n := \left\{ \gamma_{x}(t) \co x \in  p (\supp \mu )\cap \delta_n,\; t\in \left[0,1\right]  \right\}.
\]
This is a closed subset of $M$.

First let us show that
    \begin{equation}\label{rn1}
    p_{\ast}\mu (R_n) \longrightarrow 0
\end{equation}
Denote by $\chi_n$ the characteristic function of $R_n$. The sequence of functions $\chi_n$ converges pointwise to the characteristic function of $\gamma_{\gamma}(\left[0,1\right])$, so
    \[p_{\ast}\mu (R_n) \longrightarrow p_{\ast}\mu \left(\gamma_{\gamma}(\left[0,1\right])\right).
\]
Now the latter cannot be positive unless the geodesic $\gamma_{\gamma}$ is closed, for otherwise, since $\mu$ is invariant by the geodesic flow, the total mass of $\gamma_{\gamma}$ would be infinite, contradicting the fact that $\mu$ is a probability measure.
Assume $\gamma_{\gamma}$ is closed. It is two-sided because we are assuming $M$ to be orientable for the time being. Since $\supp \mu$ is closed,  $x_{\gamma}$ is in $\supp \mu$. Since $\supp \mu$ is invariant by the geodesic flow $\gamma_{\gamma}$ is contained in $\supp \mu$. Therefore by Mather's Graph Theorem
$\gamma_{\gamma}$ and $\gamma$ do not intersect. Thus by Proposition \ref{two-sided},   $\gamma$ is asymptotic to a closed geodesic, hence  cannot be  in the support of a minimizing measure by Lemma \ref{asymptote_fermee}. This proves Equation (\ref{rn1}). Besides, since  $\gamma$ is in the support of $\mu$,
\begin{equation}\label{rn2}
    p_{\ast}\mu (R_n) > 0.
\end{equation}
Next we evaluate $\inter(h,\left[\gamma_{n}\right])$ and find it equals $p_{\ast}\mu (R_n)$, which combines with the previous paragraph to prove the Proposition.

First note that by the Ergodic Decomposition Theorem (\cite{Mane}, Theorem 6.4 p. 170)
    \begin{equation}\label{rn3}
    p_{\ast}\mu (R_n)= \int_{M}\left\{\int \chi_n d\mu_x \right\}dp_\ast\mu (x)
\end{equation}
where, for $p_\ast \mu$-almost every $x$ in $M$
    \begin{eqnarray*}
    \int \chi_n d\mu_x & = & \lim_{T \rightarrow +\infty} \frac{1}{T}\int^{T}_{0} \chi_n (\gamma_x(t))dt \\
    &=& \lim_{T \rightarrow +\infty} \frac{1}{T} \sharp \left\{ t \in \left[0,T \right]\co  \gamma_x(t) \in \delta_n \right\}
\end{eqnarray*}
by the definition of $R_n$, denoting $\sharp$ the cardinal of a set.

For $x$ in $p(\supp \mu) \cap \delta_n$, let $\gamma_{x,T}$ be a closed curve obtained by closing up $\gamma_{x}(\left[0,T\right])$
with a geodesic segment $\delta_{x,T}$ of length $\leq \mbox{diam}M$. By Birkhoff's Ergodic Theorem, for $p_\ast \mu$-almost every $x$, for any closed one-form $\omega$ on $M$,
    \[\int \omega d\mu_x =  \lim_{T\rightarrow +\infty} \frac{1}{T}\int^{T}_{0} \omega_{\gamma_x(t)} (\dot\gamma_x(t))dt
    = \lim_{T\rightarrow +\infty} \frac{1}{T}<\left[\omega\right],\left[ \gamma_{x,T} \right]>.
\]
Thus, for $p_\ast \mu$-almost every $x$,
    \[\left[\mu_x \right]= \lim_{T\rightarrow +\infty} \frac{1}{T} \left[ \gamma_{x,T} \right].
\]
Since the dimension of $H_1(M,\R)$ is finite, the bilinear form $\inter (.,.)$ is continuous so for $p_\ast\mu$-almost every $x$,
    \[
    \inter (\left[\mu_x \right],\left[\gamma_n \right])=
    \lim_{T\rightarrow +\infty} \frac{1}{T}\inter (\left[ \gamma_{x,T}\right],\left[\gamma_n \right]).
\]
Observe that since both $\gamma$ and $\gamma_x$ are in the support of $\mu$, by the Graph Theorem they cannot intersect transversally. So the  transverse intersections of $\gamma_{x,T}$ and $\gamma_{n}$, if any, occur along $\delta_n$ or $\delta_{x,T}$. Note that for fixed $n$ the number $n_{x,T}$ of intersections (counted with sign) of $\delta_{x,T}$ with $\gamma$ is bounded independantly of $T$ because the length of $\delta_{x,T}$ is bounded independantly of $T$.

Furthermore, by the Graph Theorem, all  intersections of $\gamma_x (\left[0,T\right])$  with $\delta_n$ have the same sign if $\delta_n$ is small enough. This is where we need the orientability assumption.

By smoothing the corners one can make $\gamma_{x,T}$ and $\gamma_{n}$ of class $C^{1}$ without modifying their transverse intersections. The curve thus obtained are transverse unless $\gamma=\gamma_x$. In the latter case one moves $\gamma_x$ slightly away from $\gamma$ without modifying the transverse intersections of $\gamma_{x,T}$ and $\gamma_{n}$. Since all  intersections of $\gamma_x (\left[0,T\right])$  with $\delta_n$ have the same sign, we get
    \begin{equation}\label{rn4}
    \inter (\left[ \gamma_{x,T}\right],\left[\gamma_n \right])
    = \sharp \left\{ t \in \left[0,T \right]\co  \gamma_x(t) \in \delta_n \right\}+n_{x,T}
\end{equation}
whence, since $n_{x,T}$  is bounded independantly of $T$
    \[\inter (\left[\mu_x \right],\left[\gamma_n \right])=
    \lim_{T\rightarrow +\infty} \frac{1}{T}\sharp \left\{ t \in \left[0,T \right]\co  \gamma_x(t) \in \delta_n \right\}
\]
so, using Equation (\ref{rn3}),
    \[\inter(h,\left[\gamma_{n}\right])=p_{\ast}\mu (R_n)
\]
which finishes the proof of the orientable case of the Proposition.

\medskip

Assume now that $M$ is not orientable.
Let $\mu$ be a minimizing measure such that $\left[\mu\right]=rh$ with $h \in \Lambda$ and $r \in \R$. Let $\nu$ be an $I_{\ast}$-invariant  minimizing measure given by Lemma \ref{relevee_minimisante}. Let $c_1,\ldots c_b$ be an integer basis of $H^{1}(M,\R)$ and let $\omega_1,\ldots \omega_b$ be closed one-forms such that $\left[\omega_i\right]=c_i,\  i=1\ldots b$. Then $\int \omega_i d \mu \in r\Z$ for $i=1\ldots b$.

Let $\tilde{\omega}_1,\ldots \tilde{\omega}_b$ be the lifts of $\omega_1,\ldots \omega_b$ to $M_o$. They are integer one-forms and
$\left[\tilde{\omega}_1\right],\ldots \left[\tilde{\omega}_b\right]$ is a basis of $E_1 = \left\{c \in H^{1}(M_o,\R) \co I^{\ast}c=c \right\}$.
Besides,
\[
\int \tilde{\omega}_i d \nu = \int \omega_i d \mu  \;  \in r\Z,  \hskip7pt i=1, \ldots, b.
\]

Let us take an integer basis  $c_{b+1},\ldots c_{2b}$ of
$$
E_{-1} = \left\{c \in H^{1}(M_o,\R) \co I^{\ast}c=-c \right\}
$$
 and closed one-forms $\tilde{\omega}_1,\ldots \tilde{\omega}_b$ such that $\left[\tilde{\omega}_i\right]=c_i$ for $i=b+1,\ldots 2b$. Since $\nu$ is $I_{\ast}$-invariant we have
\[\int \tilde{\omega}_i d \nu = 0, \hskip7pt i=b+1, \ldots, 2b.
\]
Let $x_1,\ldots x_{2b}$ be  the coordinates of $\left[\nu \right]$ in the basis of $H_{1}(M_o,\R)$ dual to the integer basis $\left[\tilde{\omega}_1\right],\ldots \left[\tilde{\omega}_{2b}\right]$ of $H^{1}(M_o,\Z)$. We have just seen that $x_1,\ldots x_{2b}$ are all in $r\Z$, so $\left[\nu \right]$ is rational. Thus, using the orientable case of the Proposition, we conclude that $\nu$, hence $\mu$, is supported on periodic orbits.
\qed

\section{Proofs of the main theorems}
\subsection{Local results - Orientable case}
Let $h$ be a rational homology class of a surface, orientable or not. Then by Proposition \ref{rational} any $h$-minimizing measure is supported on periodic orbits. Call $\mathcal{P}_h$ the union of the projections on $M$ of the supports of all $h$-minimizing measures. By Mather's Graph Theorem $\mathcal{P}_h$ is a union of pairwise disjoint closed geodesics. Denote by $\mathcal{VP}_h$ the vector subspace of $H_{1}(M,\R)$ generated by all homology classes of geodesics contained in $\mathcal{P}_h$. Note that the convex hull of all homology classes of curves in $\mathcal{P}_h$ is contained in a flat of the unit ball containing $h$ in its interior.

The following theorem proved in \cite{gafa} describes the local geometry of the unit ball of the stable norm near a rational homology class in the orientable case.
\begin{theorem}\cite{gafa}\label{local}
Let $M$ be an orientable closed surface endowed with a Riemannian metric.
Let $h_0$ be a  rational point of $\s$. For all $h \in \mathcal{VP}_{h_0}^\perp$, there exists $s(h_0,h)>0$ such that the subset of the unit sphere $\s$
$$
 \left\{\frac{ h_0+ sh}{|| h_0+sh||} \co s \in \left[0,s(h_0,h)\right] \right\}
$$
is a straight segment.
\end{theorem}
\proof
For any $n\in \N^{\ast}$ let us denote
    \[ h_n := \frac{ h_0+ \frac{1}{n}h}{|| h_0+\frac{1}{n}h||}.
\]
Let
\begin{itemize}
    \item
    $\mu_n$ be an $h_n$-minimizing measure
    \item
    $\mu_0$ be a limit point, in the weak-$\ast$ topology, of the sequence $\mu_n$.
\end{itemize}
Then $\mu_0$ is an $h_0$-minimizing measure. By Proposition \ref{rational}, $\mu_0$ is supported on periodic orbits $\gamma_i, i\in I$ where $I$ is some set. Note that for all $i \in I$ the class $[\gamma_i]$ belongs to $\mathcal{VP}_{h_0}$. For each $i\in I$ let $V_i$ be the neighborhood of $(\gamma_i,\dot\gamma_i)$ given by Proposition \ref{two-sided}. Let $V$ be the union over $i\in I$ of the $V_i$. First let us prove that $V \cap \spt (\mu_n )$ is $\phi_t$-invariant and consists of periodic orbits homotopic to some or all of the $\gamma_i$.
Indeed by Proposition \ref{two-sided} a minimizing geodesic  that enters $V$ is either
\begin{itemize}
    \item asymptotic to one of the $\gamma_i$, which is ruled out by Lemma \ref{asymptote_fermee}
  \item homotopic to one of the $\gamma_i$
  \item or cuts one of the  $\gamma_i$ with constant sign, which is ruled out by hypothesis.
\end{itemize}
Suppose $\mu_n(V)\neq 0$. Set, for any measurable subset $A$ of $T^{1}M$
\begin{eqnarray*}
\alpha_n (A) & :=&  \frac{\mu_n(A\cap V)}{\mu_n(V)} \\
\beta_n (A) & :=&  \frac{\mu_n(A\setminus V)}{\mu_n(T^{1}M\setminus V)}\\
\lambda_n &:=& \mu_n(V).
\end{eqnarray*}
Then $\alpha_n$ and $\beta_n$ are two probability measures on $T^{1}M$. They are invariant by the geodesic flow because $V\cap \spt (\mu_n )$, as well as its complement in $ \spt (\mu_n )$, is $\phi_t$-invariant. In case $\mu_n(V)= 0$, set $\alpha_n (A)  :=  \mu_n$ and $\lambda_n := 1$. Since the support of $\alpha_n$ consists of periodic orbits homotopic to some or all of the $\gamma_i$, the homology class of $\alpha_n$ is contained in the convex hull of the
$\left[\gamma_i\right]/l_g(\gamma_i)$. Note that since the support of $\mu_0$ consists of all of the $\gamma_i$, the homology class of
$\mu_0$ is contained in the relative interior of the convex hull of the $\left[\gamma_i\right]/l_g(\gamma_i)$.

We have
    \[\mu_n = \lambda_n \alpha_n + (1-\lambda_n)\beta_n
\]
 and $\lambda_n$ tends to one as $n$ tends to infinity, so the homology class of $\alpha_n$ tends to $h_0$. Therefore, when $n$ is large enough, the homology class of
$\alpha_n$ is contained in the relative interior of the convex hull of the $\left[\gamma_i\right]/l_g(\gamma_i)$.  Thus any supporting cohomology class $c$  to $\s$ at $\left[\alpha_n\right]$, i.e.  such that $<c,\left[\alpha_n\right]> =1$ and  $<c,h> \leq 1$ for all $h \in \mathcal{B}_1$, is also a supporting cohomology class   to $\s$ at $h_0$. In other words, any flat of $\s$ that contains $\left[\alpha_n\right]$ also contains $h_0$.

Let $c$ be a supporting cohomology class to $\s$ at $h_N$. We have \linebreak $<c,h_N>=1$ and $|<c,h>| \leq 1 \  \forall h \in \s$. Therefore
    \[\lambda_N <c,\left[\alpha_N\right]> + (1-\lambda_N)<c, \left[\beta_N \right]>=1.
\]
Since $<c,\left[\alpha_N\right]> \leq 1,\ <c,\left[\beta_N \right]> \leq 1, \ \lambda_N \in \left[0,1 \right]$, this implies
$$
<c,\left[\alpha_N\right]>=<c,h_N>=1
$$
 that is, $\left[\alpha_N\right]$ and $h_N$ are in the same flat of $\s$, whence $h_0$ and $h_N$ are in the same flat of $\s$.
\qed

\bigskip

Recall from \cite{gafa} the orientable analogue of the first part of Theorem {\bf A}:

\begin{corollary}
Assume $M$ is a closed orientable surface  endowed with a Riemannian metric. Then every rational homology class contained in $\s$ lies in a flat
of $\s$ of dimension at least $b_1(M)/2-1$.
\end{corollary}
\proof Let $h_0$ be a rational point of $\mathcal{S}_1(M,g)$. Set
\[ p:= \dim \mathcal{VP}_{h_0}
\]
and assume $p<b_1(M)/2$. Choose curves
 $\gamma_i$ in $\mathcal{P}_{h_0}$ for $i=1, \ldots, p$, such that $\{\left[\gamma_i\right] \mid i=1,\ldots,p\}$ generate $\mathcal{VP}_{h}$. Since $p<b_1(M)/2$, there exists $h \in H_1(M,\R)$ such that
 \[
 h \notin \mathcal{VP}_{h_0} \mbox{ and } \inter(h,\left[\gamma_i\right])=0 \  \forall  i=1,\ldots,p.
 \]
By Theorem \ref{local} there exists $s>0$ such that
$$
F_1 := \left\{\frac{h_0+ sh}{|| h_0+sh||} \co s \in \left[0,s(h_0,h)\right] \right\}
$$
is a straight segment contained in $\mathcal{S}_1 (M, g)$.
On the other hand, the convex hull of $\left[\gamma_i \right]$ for $i=1,\ldots p+q$ is contained in a flat $F_0$ of $\mathcal{S}_1 (M,g)$ of dimension $p$ that  contains $h_0$ in its interior. From Lemma \ref{extension_plat} we deduce that there exists a flat containing $F_0$ and
$F_1$. The dimension of  said flat is greater than $p = \dim F_0$ because $h \notin \mathcal{VP}_{h_0}$.
\qed

\subsection{Local results - Non-orientable case}

In this section, we assume $M$ is a closed non-orientable surface and prove Theorems {\bf A} and {\bf B}. The proofs combine basic facts about the orientation cover of a non orientable surface and Theorem \ref{local}.

\bigskip

\begin{proposition}\label{projection_boules}
Assume $M$ is a closed non-orientable surface  endowed with a Riemannian metric $g$ and $\pi: (M_o,\tilde{g}) \rightarrow (M,g)$ its orientation cover. Then
$$
\pi_\ast \mathcal{B}_1(M_o,\tilde{g})=\mathcal{B}_1(M,g).
$$
and furthermore the  vector space $E_1$ endowed with the restriction of the stable norm of $(M_o,\tilde{g})$  is isometric to $H_1(M,\R)$ endowed with the stable norm of $(M,g)$.
\end{proposition}

\proof

Let $\mu_o$ be an element of $\mathcal{M}_{\tilde{g}}$. Then $\pi_\ast \mu_o$ is an element of  $\mathcal{M}_{g}$. So $\pi_\ast \mathcal{B}_1(M_o,\tilde{g})\subset \mathcal{B}_1(M,g)$.
Conversely, let $\mu$ be a minimizing measure of $M$.  Let $\nu  \in \mathcal{M}_{\tilde{g}}$ be given by Lemma \ref{relevee_minimisante}.
We have $\pi_{\ast}(\nu)= \mu$ so $\pi_{\ast}(\left[\nu\right])= \left[\mu \right]$. Thus $\pi_{\ast}$ restricted to  $E_1 \cap \mathcal{S}_1(M_o,\tilde{g})$ is surjective onto  $\mathcal{S}_1(M,g)$. Since $\pi_{\ast}$ is linear, it must then be surjective from $E_1 \cap \mathcal{B}_1(M_o,\tilde{g})$  onto  $\mathcal{B}_1(M,g)$. Besides, since the dimensions of $E_1$ and
$H_1(M,\R)$ are equal, $\pi_{\ast}$ restricted to  $E_1$ must be injective. So $\pi_{\ast}$ restricted to  $E_1$ is a linear isomorphism sending $E_1 \cap \mathcal{B}_1(M_o,\tilde{g})$ to $\mathcal{B}_1(M,g)$.
\qed

\bigskip

The purpose of the next Proposition is to evaluate the maximal dimension of a flat containing a rational class $h$ (not necessarily as an
interior point), depending on the topological properties of $h$-minimizing curves. Recall that a simple closed curve $\gamma$ of $M$ is said
{\it of type I} (resp. {\it of type II}) if its inverse image $\pi^{-1}(\gamma)$ consists of either one curve or two homologous curves (resp.
two non-homologous curves). Let $h$ be a rational point of $\mathcal{S}_1(M,g)$. Partition $\mathcal{P}_h$ in two subsets $\mathcal{P}^{1}_{h}$
and $\mathcal{P}^{2}_{h}$, the former  consisting only of curves of type I and the latter only of curves of type II. Let $\mathcal{VP}^{2}_{h}$
be the vector subspace of $H_{1}(M,\R)$ generated by all homology classes of geodesics contained in $\mathcal{P}^{2}_{h}$. Let
$\mathcal{VP}^{1}_{h}$ be such that $\mathcal{VP}^{1}_{h}$ is generated by homology classes of curves of type I and
\[ \mathcal{VP}^{2}_{h} \oplus \mathcal{VP}^{1}_{h}= \mathcal{VP}_{h}.
\]

\begin{proposition} \label{nonor}
Let $M$ be a closed non-orientable surface and let $h_0$ be a rational point of $\mathcal{S}_1(M,g)$. Set
\[ p:= \dim \mathcal{VP}^{1}_{h_0} \mbox{ and } q:= \dim \mathcal{VP}^{2}_{h_0}
\]
and assume $p+2q<b_1(M)$. Then there exists a flat of $\mathcal{B}_1(M,g)$ containing $h_0$, of dimension $>p+q$.
\end{proposition}
\proof
Choose curves
\begin{itemize}
\item $\gamma_i$ in $\mathcal{P}^{1}_{h_0}$ for $i=1, \ldots, p$, such that $\{\left[\gamma_i\right] \mid i=1,\ldots,p\}$ generate $\mathcal{VP}^{1}_{h}$
\item $\gamma_i$ in $\mathcal{P}^{2}_{h_0}$ for $i=p+1, \ldots, p+q$, such that $\{\left[\gamma_i\right] \mid  i=p+1,\ldots,p+q\}$ generate $\mathcal{VP}^{2}_{h}$.
\end{itemize}
Denote by $\mu_i$ the $\phi_t$-invariant probability measure supported on $\gamma_i$ for $i=1,\ldots p+q$. Let $c \in H^{1}(M,\R)$ be such that $h_0$ is $c$-minimizing. Then each $\mu_i$, and each convex combination thereof, is also $c$-minimizing.  Let $\lambda_i \in \left]0,1\right[$, $i=1,\ldots p+q$ be such that $\sum_i \lambda_i = 1$ and $\sum_i \lambda_i \left[\mu_i\right]=h_0$.

\medskip

If $i \in \left\{1,\ldots p\right\}$, choose a closed geodesic $\tilde{\gamma}_i$ in $M_o$ such that
\begin{itemize}
    \item $\pi(\tilde{\gamma}_i)=\gamma_i$
    \item $\left[\tilde{\gamma}_i\right]$ lies in the eigenspace $E_1$ for the involution $I$.
\end{itemize}

 If $i \in \left\{p+1,\ldots p+q\right\}$, choose two closed geodesics $\tilde{\gamma}_i$ and $\tilde{\gamma}_{i+q}$ in $M_o$ such that
\begin{itemize}
    \item $\pi(\tilde{\gamma}_i)=\gamma_i$
    \item $I(\tilde{\gamma}_i)= \tilde{\gamma}_{i+q}$
    \item $\left[\tilde{\gamma}_i\right]\neq  \left[\tilde{\gamma}_{i+q}\right]$.
\end{itemize}

\medskip

Define
\begin{itemize}
    \item $\tilde{\mu}_i$ the $\tilde{\phi}_t$-invariant probability measure supported on $\tilde{\gamma}_i$ for $i=1,\ldots, p+2q$
    \item $\tilde{\lambda}_i := \lambda_i$ if $i=1,\ldots, p$
    \item $\tilde{\lambda}_i := \lambda_i /2$ if $i=p+1,\ldots, p+q$
    \item $\tilde{\lambda}_{i+q} := \lambda_i /2$ if $i=p+q+1,\ldots, p+2q$
    \item $\tilde{\mu} := \sum^{p+2q}_{i=1} \tilde{\lambda}_i \tilde{\mu}_i$
    \item $\tilde{h}_0 = \left[\tilde\mu\right]$.
\end{itemize}

We have
\begin{itemize}
    \item $I_{\ast}(\tilde{\mu})=\tilde{\mu}$ whence $I_{\ast}(\tilde{h}_0)=\tilde{h}_0$
    \item $\pi_{\ast}(\tilde{\mu})=\mu$ whence $\pi_{\ast}(\tilde{h}_0)=h_0$
    \item $\tilde{\mu}$ is $\pi^{\ast}(c)$-minimizing
    \item the vector space generated by $\left[\tilde{\gamma}_i\right]$ for $i=1,\ldots p+2q$ equals $\mathcal{VP}_{\tilde{h}_0}$.
\end{itemize}
The last equality stands because of Proposition \ref{projection_boules}. To clear up the notation, call $V$ the vector subspace of
$H_{1}(M_o,\R)$ generated by the integer classes $\left[\tilde{\gamma}_i\right]$  for $i=1,\ldots, p+2q$. Note that $I(V)=V$, so $V=V_1 \oplus
V_{-1}$ where $V_i = E_i \cap V$, $i=\pm 1$. Also, $V = \mathcal{VP}_{\tilde{h}_0}$. We have
\begin{eqnarray*}
V_1 &=& \mbox{Vect}\left(\left\{ \left[ \tilde{\gamma}_i \right] \co i=1,\ldots p \right\}\cup
    \left\{ \left[ \tilde{\gamma}_i \right] + \left[ \tilde{\gamma}_{i+q} \right] \co i=p+1,\ldots p+q \right\}\right)\\
V_{-1} &=& \mbox{Vect}\left(\left\{ \left[ \tilde{\gamma}_i \right] - \left[ \tilde{\gamma}_{i+q} \right] \co i=p+1,\ldots p+q \right\}\right)
\end{eqnarray*}

We would like to use Theorem \ref{local} with $\tilde{h}_0$ playing the part of $h_0$ and some $h$ in $E_1 \cap V^\perp$ but not in $V$.
Observe that
    \begin{eqnarray*}
     \dim V^{\perp}_{-1} & = &  b_1(M_o) - q \mbox{ so }\\
     \dim V^{\perp}_{-1}\cap E_1 & \geq &  b_1(M_o) - q + b_1(M) -b_1(M_0)\\
     &=&  b_1(M) - q > p+q = \dim V_1
  \end{eqnarray*}
since we assume $b_1(M)  > p+2q$. So there exists $h \in V^{\perp}_{-1}\cap E_1$ such that $h \notin V_1$. Since $V_1 \subset E_1$, we have $E_1 =  E^{\perp}_{1} \subset  V^{\perp}_{1}$ thus $h \in V^{\perp}_{1}$ and
\[ h \in V^{\perp}_{1}\cap V^{\perp}_{-1}= \left(V_1 \oplus V_{-1}\right)^{\bot}=V^{\bot}.
\]
So $h \in E_1 \cap V^\perp$. Furthermore $h \notin V$ since $h \in E_1$ and $h\notin V_1 = E_1 \cap V$.

\medskip

By Theorem \ref{local} there exists $s>0$ such that
$$
F_1 := \left\{\frac{ \tilde{h}_0+ sh}{|| \tilde{h}_0+sh||} \co s \in \left[0,s(h_0,h)\right] \right\}
$$
is a straight segment contained in $\mathcal{S}_1 (M_o, \tilde{g})\cap E_1 \cong \mathcal{S}_1 (M,g)$.
On the other hand, the convex hull of $\left[\gamma_i \right]$ for $i=1,\ldots p+q$ is contained in a flat $F_0$ of $\mathcal{S}_1 (M,g)$ that  contains $h_0$ in its interior. From Lemma \ref{extension_plat} we deduce that there exists a flat containing $F_0$ and
$F_1$. The dimension of such a flat is greater than $p+q = \dim F_0$ because $h \notin V$.

\qed

\bigskip

Taking a rational $h$ in the proof of Proposition \ref{nonor}, we deduce  the first part of Theorem {\bf A}.

\begin{corollary}
Assume $M$ is a closed non-orientable surface  endowed with a Riemannian metric. Then every connected minimizing cycle is a component of a minimizing cycle whose homology class lies in a flat of $\s$ of dimension at least $[(b_1(M)+1)/2]-1$.
\end{corollary}

\bigskip

Let $\Gamma$ be a minimizing cycle whose connected components  are not pairwise proportional in homology. Then its connected components form a
system of disjoint non pairwise homologically proportional simple closed curves. We have the following result.

\begin{proposition}\label{curves}
A system of disjoint non pairwise homologically proportional simple closed curves has its cardinality bounded from above by $2b_1(M)-1$.
\end{proposition}

\proof The argument is classical. Let $\alpha_1,\ldots,\alpha_{p+q}$ be a maximal system of disjoint pairwise non homologically proportional
simple closed curves of $M$. Suppose that $p$ is the number of one-sided curves of this system and $q$ the number of two-sided curves. By
cutting $M$ along these simple closed curves we obtain an union of $b_1(M)-1$ pair of pants. So we must have $3(b_1(M)-1)=p+2q$. This implies
$2(p+q)=3.b_1(M)+p-3$. As $p\leq b_1(M)+1$, the assertion follows. \qed

\medskip

\noindent So $\Gamma$ has at most $2b_1(M)-1$ components. The second part of theorem {\bf A} is proved.

\bigskip

By specializing Proposition \ref{nonor} to the case when $q=0$, we now deduce  Theorem {\bf B} which describes the local geometry of the unit ball of the stable norm near a rational homology class for which the connected components of minimizing cycles are curves of type I.

\bigskip

\begin{corollary}
Assume $M$ is a closed non-orientable surface endowed with a Riemannian metric. Let $h_0$ be an integer homology class all of whose minimizing cycles consist of curves of type I. Then for all $h \in H_1(M,\R)$, there exists $s(h_0,h)>0$ such that the subset of the unit sphere $\s$
$$
 \left\{\frac{ h_0+ sh}{|| h_0+s h||} \co s \in \left[0,s(h_0,h)\right] \right\}
$$
is a straight segment.
\end{corollary}

\bigskip

\subsection{Global result - Proof of Theorem C}

We first prove the following theorem.

\begin{theorem} \label{pol}

Let $M$ be a closed surface (orientable or not) endowed with a Riemannian metric $g$ and $c_1,\ldots,c_l$ a family of disjoint smooth, simple, closed curves whose homology classes are not pairwise proportional (that is $[c_i] \notin \R [c_j]$ for $i\neq j$).

For all sequence $\{r_i\}_{i=1}^l$ of positive real numbers, there exists a smooth metric $g^\ast$ conformal to $g$ such that the intersection of $\s(g^\ast)$ with the subspace spanned by the curves $[c_1],\ldots,[c_l]$ coincides with the polyhedron  $$\text{Conv}_s\left(\frac{[c_1]}{r_1}, \ldots,\frac{[c_l]}{r_l}\right)$$ where $\text{Conv}_s$ denote the convex hull of the symmetrized of a set.
\end{theorem}

\proof Let $P$ denote the polyhedron generated as the convex hull
$$
\text{Conv}_s\left(\frac{[c_1]}{r_1}, \ldots,\frac{[c_l]}{r_l}\right).
$$
We can suppose that each curve $[c_i]$ corresponds to an exposed point of the polyhedron (if not we can discard this curve and the polyhedron $P$ remains unchanged).

\begin{lemma} \label{lemmapol1}
There exists a smooth metric $\bar{g}$ conformal to $g$ and an open neighborhood $V_i$ of each $c_i$ such that $c_i$ is the unique closed $\bar{g}$-geodesic of $V_i$ and $l_{\bar{g}}(c_i)=r_i$.
\end{lemma}

\proof[Proof of Lemma \ref{lemmapol1}.]
For $\epsilon$ small enough, the $\epsilon$-tubular neighborhoods $U_\epsilon(c_i)$ are pairwise disjoints and the $g$-orthogonal projections $p_i : U_\epsilon(c_i)\rightarrow c_i$ are well defined. For each $x \in U_\epsilon(c_i)$, there exists two $g$-unitary vectors $\pm v(x) \in T_x M$  orthogonal to the fiber $p_i^{-1}(p_i(x))$. The function $\alpha_i :U_i \rightarrow ]0,\infty[$ given by the formula
$$
\alpha_i(x)=g(Dp_i(x)v(x),Dp_i(x)v(x))
$$
is smooth and such that $\alpha_i \circ c_i=1$. We define a new metric $g'$ conformal to $g$  by $\alpha_i g$ on $U_\epsilon(c_i)$ and by extending the local conformal factors $\{\alpha_i\}_{i=1}^k$ into a smooth positive function $\alpha$ on M. We claim that the projections $p_i : U_\epsilon(c_i) \rightarrow c_i$ do not increase the lengths with respect to $g'$. Indeed, take $x \in U_i$ and $w \in T_x U_i$. Write $w = \lambda v(x) +\mu v'$, with $v' \in T_x p_i^{-1}(p_i(x))$. Note that $v$ and $v'$ are $g$-orthogonal ; since $g'$ is conformal to $g$, $v$ and $v'$ are $g'$-orthogonal. Orthogonal projections do not increase distances, so
    \[\lambda^{2}g'(v(x),v(x))= g'(\lambda v(x),\lambda v(x)) \leq g'(w,w).
\]
Now
    \[g'(Dp_i(x) w,Dp_i(x) w)=g'(Dp_i(x) \lambda v,Dp_i(x) \lambda v)
\]
because $p_i$ is the orthogonal projection to $c_i$, whence

\begin{eqnarray*}
g'(Dp_i(x) w,Dp_i(x) w)& = & \lambda^{2}\alpha_i(x)\\
    & = & \lambda^{2}\alpha_i(x)g(v,v)\\
    & = & \lambda^{2}g'(v,v) \leq g'(w,w) \\
\end{eqnarray*}
which proves the claim.

 Choose a function $f \in {\mathcal{C}}^\infty(M)$ null on $\cup_{i=1}^k c_i$,  positive elsewhere and such that $\Delta(f)<-\text{Scal}_{g'}$ where $\text{Scal}_{g'}$ denote the scalar curvature of $g'$. Let $g'':=\exp(f)g'$. We can easily verify that each projection $p_i : U_\epsilon(c_i) \rightarrow c_i$ now strictly contracts the lengths. So $c_i$ is a $g''$-geodesic and the negativity of $\text{Scal}_{g''}=\exp(f)\Delta(f)+\text{Scal}_{g'}$ ensures that the orbit $c_i$ of the geodesic flow associated to $g''$ is hyperbolic, hence isolated  in a neighborhood $V_i$.

We extend the functions $\lambda_i=r_i / l_{g}(c_i)$ defined on each neighborhood $V_i$ into a smooth function $\lambda$ defined on the whole surface and set $\bar{g}:=\lambda^2 g''$. The lemma is proved. \qed

\medskip

For any sequence $\eps$ such that $\epsilon_i=0,\pm 1$ we denote by $\gamma(\eps)$ the multicurve $\cup_{i=1}^l \epsilon_i \cdot c_i$ of $\cup_{i=1}^l V_i$ minimizing the length in the class $\sum_{i=1}^l \epsilon_i [c_i]$.

\begin{lemma}\label{lemmapol2}
There exists a smooth metric $g^\ast$ conformal to $\bar{g}$ such that for any sequence $\eps:=\{\epsilon_i\}_{i=1}^l$ with $\epsilon_i=0,\pm 1$,
$$
||\sum_{i=1}^l \epsilon_i [c_i]||_s^{g^\ast}=l_{g^\ast}(\gamma(\eps)).
$$
\end{lemma}

\proof[Proof of Lemma \ref{lemmapol2}.] Set
$$
\delta(\eps):=l_{\bar{g}}(\gamma(\eps))-||\sum_{i=1}^l \epsilon_i [c_i]||_s^{\bar{g}}.
$$
The set $\Gamma(\eps)$ of unions of closed geodesics different from $\gamma(\eps)$ homologous to $\sum_{i=1}^l \epsilon_i [c_i]$ such that their length is bounded from above by $l_{\bar{g}}(\gamma(\eps))$ is compact. It is clear that no multicurve $\gamma$ in $\Gamma(\eps)$ can  be totally contained in $\cup_{i=1}^l V_i$. So, if $\delta(\eps)>0$, the infimum $t(\eps)$ of time $t$ such that there exists a multicurve $\gamma \in \Gamma(\eps)$ spending a time $t$ outside $\cup_{i=1}^l V_i$ is reached and not zero.

Since there is but a finite number of sequences $\eps$ we  may choose $\beta$ so big  that for all $\eps$ with $\delta(\eps)>0$
$$
\beta>\log\left[1+\frac{\delta(\eps)}{t(\eps)}\right].
$$
Now we choose a function $f' \in {\mathcal{C}}^\infty(M)$ null on $\cup_{i=1}^l c_i$,  positive elsewhere and such that $f'>\beta$ outside $\cup_{i=1}^l V_i$. For any multicurve $\gamma$ spending some time $t$ outside $\cup_{i=1}^l V_i$, we have
$$
l_{\exp(f')\bar{g}}(\gamma)>(\exp(\beta)-1)t+l_{\bar{g}}(\gamma).
$$

Now let $g^\ast=\exp(f')\bar{g}$. All the lengths except those of the $c_i$'s increase for $g^\ast$ so for all $\eps$ such that $\delta(\eps)>0$ and for all multicurve $\gamma$ in the class $\sum_{i=1}^l \epsilon_i [c_i]$,
$$
l_{g^\ast}(\gamma)\geq l_{g^\ast}(\gamma(\eps)).
$$
\qed

\medskip

Recall that $l_{g^\ast}(c_i)=r_i$ as the lengths of the $c_i$'s do not increase. Thus each exposed point of $P$ belongs to the unit sphere $\s(g^\ast)$ of the stable norm. Furthermore by lemma \ref{lemmapol2} each face of $P$ contains an interior point that belongs to $\s(g^\ast)$. This proves the theorem by convexity of the unit sphere of the stable norm. \qed

\bigskip

\begin{figure}[h]
\leavevmode \SetLabels
\L(.49*.4) $M$\\
\endSetLabels
\begin{center}
\AffixLabels{\centerline{\epsfig{file =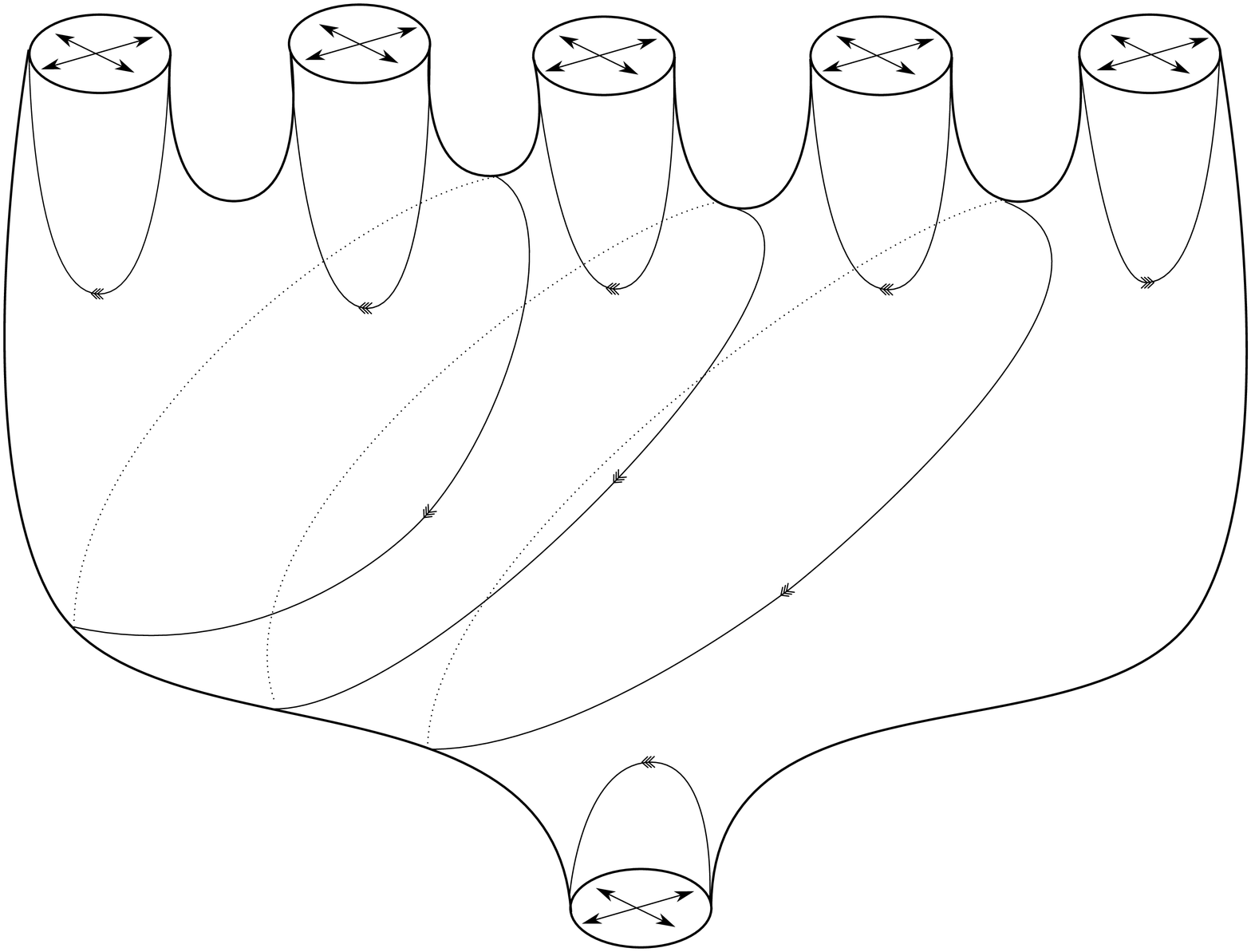,width=7cm,angle=0}}}
\end{center}
\caption{}
\label{famille1}
\end{figure}

\begin{corollary}
Let $M$ be a closed and orientable surface endowed with a Riemannian metric $g$. For each isotropic subspace $L$ of $H_1(M,\R)$ (with respect to $\inter$), there exists a metric $g^\ast$ conformal to $g$ such that the restriction of $\s(g^\ast)$ to $L$ is a polyhedra with rational vertices.
\end{corollary}

\proof There exists a system of disjoint, smooth, simple and non-pairwise homotopic closed curves $c_1,\ldots,c_l$ that span $L$. We apply Theorem \ref{pol} to obtain the claim. \qed

\medskip

\noindent {\bf Remark.} It is a classical result that such a system has cardinality at most $(3/2)b_1(M)-3$ (same argument  as in the proof of Proposition \ref{curves}), thus this bounds the number of vertices of the polyhedra obtained that way by $3b_1(M)-6$.

\bigskip

We now deduce, as a corollary of Theorem \ref{pol}, Theorem {\bf C} as stated in the introduction :

\begin{corollary}
Let $M$ be a closed and non-orientable surface endowed with a Riemannian metric $g$. There exists a metric $g^\ast$ conformal to $g$ such that $\s(g^\ast)$ is a polyhedra with rational vertices.
\end{corollary}

\proof  There exists a system $c_1,\ldots,c_l$ of smooth, simple, closed curves such that $c_i \cap c_j=\emptyset$, $[c_i] \notin \R [c_j]$ for $i\neq j$ and $H_1(M,\R)=\text{Vect}([c_1],\ldots,[c_l])$ (see figure \ref{famille1} for an example of such a system with $l=2b_1(M)-1$). We apply  Theorem \ref{pol} to obtain the claim. \qed

\medskip

\noindent {\bf Remark.} Such a system has cardinality at most $2b_1(M)-1$ (proposition \ref{curves}), thus this bounds the number of vertices of the polyhedra obtained that way by $4b_1(M)-2$.

\bigskip

\noindent {\bf Acknowledgements.} The authors are grateful to Ivan Babenko for several helpful conversations.

{\small

\bigskip

\noindent Florent Balacheff\\
Section de Math\'ematiques, Universit\'e de Gen\`eve, Suisse\\
e-mail : florent.balacheff@math.unige.ch

\medskip

\noindent Daniel Massart\\
D\'epartement de Math\'ematiques, Universit\'e Montpellier II, France\\
e-mail : massart@math.univ-montp2.fr
}

\end{document}